\documentclass[english, 12pt]{article}
\usepackage{amsfonts}
\usepackage{amsmath}
\usepackage{amssymb}
\usepackage{amsthm}
\usepackage{amscd}
\usepackage{amscd}
\usepackage{amsthm}

\usepackage[
  backend = bibtex8 
  ,style = numeric 
  ,natbib = false 
  ]{biblatex}
\usepackage{amsfonts}
\usepackage{amsmath}
\usepackage{color}
\usepackage{mathrsfs}

\newtheorem{defi}{Definition}[section]
\newtheorem{theo}{Theorem}[section]
\newtheorem{prop}{Proposition}[section]
\newtheorem{lem}{Lemma}[section]

\newtheorem{rk}{Remark}[section]

\newtheorem{coro}{Corollary}[section]
\numberwithin{equation}{section}

\def\calE{{\cal{E}}}

\def\calE{\mathcal{E}}

\def\R{\mathbb{R}}
\def\N{\mathbb{N}}

\def\alp{\alpha}
\def\Om{\Omega}

\def\P{{\bf P}}

\DeclareMathOperator{\ran}{ran}
\DeclareMathOperator{\dom}{dom}

\newcommand{\rel}{\mathrm{rel}}

\begin{document}
	\title{Feller property and convergence for semigroups of time-changed processes}
	\author{\normalsize Ali BenAmor\footnote{High school for transport and logistics, University of Sousse, Sousse, Tunisia. E-mail: {\tt ali-ben05} {\tt @yahoo.de}}
	 \ and Kazuhiro Kuwae\footnote{Department of Applied Mathematics, Faculty of Science, Fukuoka University, Fukuoka 814-0180, Japan. E-mail: {\tt kuwae@fukuoka-u.ac.jp} 
 }	
	}
	
\date{}
	
\maketitle
	
{\bf Abstract:} We give a substitute to the Feller property for semigroups of time-changed processes; under some conditions this leads to  sufficient (new) conditions for the semigroups to be Feller. Moreover, given a standard process and a sequence of measures converging vaguely to a final measure, under some assumptions,  we establish convergence of the sequence of the semigroups and the resolvents  of the corresponding  time changed-processes. Some  applications are given: convergence of solutions of evolution equations and convergence of finite time distributions, as well as weak convergence of the related processes. \\
%

\section{Introduction}
Let $X$ be  a locally compact separable metric space and $m$ a positive Radon measure on Borel subsets of $X$ with full support $X$. Let us consider a  standard process $\mathbb{X}$ with state space $X$, $\mu$ a positive Borel smooth measure on $X$ with fine support $F$ and $\bf{A}$ the PCAF of $\mu$ (via Revuz correspondence). The main objects under study in this paper are the resolvent and the semigroup of the  time changed process $\mathbb{X}$ by $\bf{A}$ (see details below) which we denote by  $\check{\mathbb{X}}$, being considered on the state space $X$. It is known that $\check{\mathbb{X}}$ considered on $X$ does not give rise to a right process, and hence neither its resolvent, which we denote by $\check{R}_\alp$, needs to be  strongly continuous, nor its semigroup, which we denote by $\check{\bf P}_t$ needs to be a ${\bf C_0}$-semigroup on the spaces of continuous bounded functions or of continuous functions  vanishing at infinity. Even worst: $\check{R}_\alp$ is not injective and $\check{\bf P}_t$ is a degenerate semigroup, i.e., there are non-vanishing functions such that $\check{\bf P}_tu = 0$ for all $t\geq 0$.\\
Much is known about $\check{\mathbb{X}}$ with state space the fine support of the measure $\mu$, i.e. on the set $F$, or its topological support, especially in the $L^2$-setting (here we refer to \cite{CF}). However, few  is known about $\check{\mathbb{X}}$ with state space $X$. For instance, very few results are known concerning  continuous dependence of both the resolvent and the semigroup with respect to variations of the measure $\mu$. Actually, if $X$ is compact and the Green function satisfies some technical conditions, Duncan \cite{Duncan} proved continuous dependence of the semigroup with respect to variation of the additive functional.\\
Our aim in these notes is twofold.  First, we consider $\check{\mathbb{X}}$ on $X$ and analyse the properties of the related resolvent and semigroup on the space of continuous functions on $X$ vanishing at infinity, $C_0(X)$, under the condition that $\mu$ is a $G$-Kato measure. We will show that the out-coming generator is a Hille--Yosida operator and hence the generated  semigroup is a ${\bf C_0}$-semigroup and the resolvent is strongly continuous on a closed subspace of $C_0(X)$ that we determine explicitly. As a by-product we obtain  a substitute for Feller property for $\check{\mathbb{X}}$; this property turns to be the classical Feller property if $F=X$. Accordingly, our result establishes a procedure on how to get Feller processes by changing the reference measure: If $F=X$ and $\mu$ is a $G$-Kato measure then $\check{\mathbb{X}}$ considered on $X$ has the Feller property. Let us stress at this stage that compared with \cite[Section 7]{Kuwae-17}, we do not assume the starting process $\mathbb{X}$ to be a Feller process and that the condition $F=X$ is still relevant theoretically and practically (see also \cite{Kuwae-24}).\\
Besides we show that $\check{\bf P}_t$ induces  a $1$-integrated Lipschitz-continuous semigroup on the space $C_0(X)$  (not to confound with $\bf{C_0}$-semigroup!). This result can be  exploited to solve the Cauchy problem  for a larger class of initial data.\\
We also give necessary and sufficient conditions for $\check{\mathbb{X}}$ to be a standard process on $X$. Here we mainly consider the normality property, because all other properties of a standard process are inherited from $\mathbb{X}$.

Our second aim in this study is to analyse into which extent  the resolvent and the semigroup of $\check{\mathbb{X}}$ depend continuously  on variations of the underlying measure $\mu$. Precisely,  let $(\mu_n)$ be a sequence  of measures converging vaguely to a measure $\mu_\infty$. Let $\bf{A}^{\mu_n}$ be the PCAF of $\mu_n$ and $\check{\mathbb{X}}^n$ be the resulting time-changed process obtained from $\mathbb{X}$, for each $n\in\N\cup\{\infty\}$. Under some specific conditions, which are fulfilled for a large class of situations, we show that the related resolvents as well as parts of the related semigroups  considered on the space $C_0(X)$ behave continuously in the strong topology. Moreover stronger conditions on the fine supports or monotone convergence of the $\mu_n$'s lead even to continuous dependence of the semigroups. As a by-product we obtain convergence of the related $1$-integrated semigroup as well as convergence and approximation for solutions of some evolution equations.

Our study is motivated, first, by the fact that the generator of $\check{\mathbb{X}}$ with state space $X$ is still a Hille--Yosida operator and hence still defines a ${\bf C}_0$-semigroup  and a strongly continuous resolvent, however on a subspace of $C_0(X)$, which enables one to solve the first order Cauchy problem in the space $C_0(X)$ with some particular initial data. Furthermore, owing to recent results due  to Arendt--Chalendar--Moletsane \cite{Arendt-23} we will establish that the generator, which is in fact a linear relation,  generates  a $1$-integrated semigroup. Thereby, we get a generic example (coming from probability) of $1$-integrated semigroups; which are important technical tool to solve the  Cauchy problem on the space of continuous functions.\\
The second motivation  lies in the role played by convergence of resolvent and semigroup for solving stationary and evolution equations as they express continuous dependence of solutions on initial data and lead also to an approximation procedure. Besides such problem is not yet intensively studied for the time-changed process with state space $X$. In this respect we mention Duncan's paper  \cite{Duncan}, in which he studied continuous dependence of the resolvent and semigroup with respect to the additive functional  on the state space $X$, however under the assumption that $X$ is compact together with assumptions on the kernel of the resolvent which imply that the measure is $G$-Kato in our sense (see \cite[(A) and Poposition 2.1]{Duncan}).\\
We also refer to \cite{BenAmor-24} where continuous dependence of spectral data with respect to monotone weak convergence of the $\mu_n$'s is studied.

\section{Preparing results}
Let $X$ be a locally compact separable metric space, $m$ a positive Radon measure on Borel subsets of $X$ with full support. 
The symbol $\mathscr{B}(X)$ for Borel $\sigma$-field over $X$ is also used as the space of Borel measurable functions on $X$. 
We designate by $\mathscr{B}_b(X)$ resp. $\mathscr{B}_b^+(X)$ the spaces of Borel bounded and positive Borel bounded functions on $X$ endowed with the topology of uniform convergence which we denote by $\|\cdot\|_\infty$.\\
Set $X_{\partial}:= X\cup\{\partial\}$ the one point compactification of $X$, $\mathscr{B}_b(X_\partial)$ the space of numerical Borel measurable bounded functions on $X_\partial$ and $\mathscr{B}^+(X_\partial)$ the space of numerical Borel measurable positive functions on $X$. We will extend any function $u\in\mathscr{B}(X)$ to $X_\partial$ by setting $u(\partial) = 0$ so that the extended function is Borel measurable on $X_\partial$ and one can identify $\mathscr{B}(X)$ with the subset of functions in $\mathscr{B}(X_\partial)$ vanishing at $\partial$. 
We denote by $\mathscr{B}_b^*(X)$ the set of bounded universally measurable functions on $X$, where the $\sigma$-field $\mathscr{B}^*(X)$ over $X$ is the universal completion of $\mathscr{B}(X)$ (see \cite[A.1]{Sharpe} for the definitions of 
universal completion of $\mathscr{B}(X)$).
\\
Let 
\[
	\mathbb{X}:=	\left( \Om,\mathscr{M}, \mathscr{M}_t, (\mathbb{P}_x)_{x\in X_{\partial}}, (\mathbb{X}_t)_{0\leq t\leq\infty}, (\theta_t)_{0\leq t\leq\infty},\zeta \right),
\]
be a standard process in the sense of \cite[Definition (9.2), p.~45]{Blumenthal-Getoor} (not necessarily symmetric) with state space $(X,\mathscr{B}(X))$ (augmented by $\partial$). Here $(\Om,\mathscr{M})$ is a measurable space, $(\mathbb{P}_x)_{x\in X_\partial}$ is a family of probability measures on $\Om$ and $\mathbb{X}_t : \Om\to X_\partial$ are such that $(\mathbb{X}_t)$ is a stochastic process on $( \Om,\mathscr{M} )$. 
\\
We first assume that $\mathbb{X}$ is transient, unless otherwise explicitly stated, and it has a dual process $\hat{\mathbb{X}}$ with the same properties.\\
For any $u\in \mathscr{B}_b(X_\partial)\cup \mathscr{B}^+(X_\partial), t\geq 0$ and $\alpha>0$ set
\begin{align*}
		{\bf P}_tu(x) : = \mathbb{E}_x\left[ u(\mathbb{X}_t) \right], \quad x\in X_\partial.
\end{align*}	
Clearly if $u\in \mathscr{B}_b(X)\cup \mathscr{B}^+(X)$  then with our convention we have
\begin{align*}
		{\bf P}_tu(x) : = \mathbb{E}_x\left[ u(\mathbb{X}_t); \mathbb{X}_t\in X \right],\quad x\in X.
\end{align*}
Moreover ${\bf P}_tu(\partial) = 0$ and our convention is consistent. Finally we quote that $({\bf P}_t)_{t>0}$ is a semigroup of contractions on $\mathscr{B}_b(X)$; it is the transition semigroup of $\mathbb{X}$.\\
Right-continuity of the sample paths imply that for any $u\in\mathscr{B}_b(X)$ the map $(t,x)\mapsto {\bf P}_t u(x)$ is bounded w.r.t. $x$ and $\mathscr{B}([0,\infty))\times\mathscr{B}(X)$-measurable, where $\mathscr{B}([0,\infty))$ is the Borel algebra on $[0,\infty)$ (see \cite[p.~41]{Blumenthal-Getoor}).  Thus by \cite[Theorem 1.5, p.~35]{Dynkin}, for any $u\in\mathscr{B}^+(X) \cup \mathscr{B}_b^*(X), \alp>0$ one can define the Bochner integral
\begin{align*}	
	R_\alpha u(x) := \mathbb{E}_x\left[\int_0^\infty e^{-\alpha t} u(\mathbb{X}_t)\,dt \right]
		= \int_0^\infty e^{-\alpha t}  {\bf P}_t u(x)\,dt,\quad x\in X.
\end{align*}
Obviously we have $R_\alp :\mathscr{B}_b(X)\to \mathscr{B}_b(X)$ and  $\alpha\|R_\alpha u\|_\infty\leq \|u\|_\infty$ for $u\in\mathscr{B}_b(X)$. Moreover, by the semigroups property we obtain the resolvent identity $R_\alp - R_\beta = (\beta - \alp)R_\alp R_\beta$. The family $(R_\alpha)_{\alp>0}$ is the resolvent of $({\bf P}_t)_{t>0}$.\\ 
A set $B\subset X_{\partial}$ is said to be a \emph{nearly Borel set} (relative to the given process $(\mathbb{X}_t)$) if for each Borel probability measure $\mu$ on $X_{\partial}$ there exist Borel subsets $B_1$ and $B_2$  of $X_{\partial}$ such that $B_1\subset B\subset B_2$ and 
\begin{align*}
\mathbb{P}_{\mu}(\mathbb{X}_t\in B_2\setminus B_1\text{ for some }t\geq0)=0.
\end{align*}
It is known that for any nearly Borel set $B$, its first hitting time $\sigma_B:=\inf\{t>0: \mathbb{X}_t\in B\}$ is an $(\mathscr{M}_t)$-stopping time (see \cite[p.~408]{CF}). Denote by $\mathscr{B}^{\,n}(X)$ the family of nearly Borel subsets of $X$, or nearly Borel measurable functions on $X$. The symbol $\mathscr{B}^{\,n}(X_{\partial})$ can be defined similarly.\\   
A set $N\subset X$ is said to be \emph{thin} if there exists a set $B\in\mathscr{B}^{\,n}(X)$ such that $N\subset B$ and 
$B^r=\emptyset$, where $B^r:=\{x\in X:\mathbb{P}_x(\sigma_B=0)=1\}$. 
A set $N\subset X$ is said to be \emph{semi-polar} if $N$ is contained in a countable union of thin sets. 
It is known that for any $B\in\mathscr{B}^{\,n}(X)$, $B\setminus B^r$ is semi-polar (\cite[Chapter II, Proposition~3.3]{Blumenthal-Getoor}) and the set $\{t\in[0,\infty):\mathbb{X}_t\in B\}$ is at most countable $\mathbb{P}_x$-a.s.~for all $x\in X$ (\cite[Chapter II, Proposition~3.4]{Blumenthal-Getoor}). 
A set $N\subset X$ is said to be \emph{thin at $x$} if there exists a set $B\in\mathscr{B}^{\,n}(X)$  such that $N\subset B$ and 
$x\not\in B^r$.
A set $O\subset X_{\partial}$ is said to be \emph{finely open} if 
$X_{\partial}\setminus O$ is thin at each $x\in O$, in other words for each $x\in O$ there exists a set $B\in\mathscr{B}^{\,n}(X_{\partial})$ such that $X\setminus O\subset B$ and 
$x\notin B^r$. Let $\mathscr{O}(X)$ ($\mathscr{O}(X_{\partial})$) be the collection of all finely open subsets of $X$ ($X_{\partial}$). Then $\mathscr{O}(X)$ ($\mathscr{O}(X_{\partial})$) is a topology on $X$ ($X_{\partial}$), which is called the \emph{fine topology} on $X$ ($X_{\partial}$). For a nearly Borel subset $F$ of $X$, $F$ is finely closed if and only if $F^r\subset F$. 
A set $N\subset X$ is called \emph{exceptional} or \emph{$m$-polar} if there exists a set $\widetilde{N}\in\mathscr{B}^{\,n}(X)$ with $N\subset \widetilde{N}$ such 
that $\mathbb{P}_m(\sigma_{\widetilde{N}}<\infty)=0$. 
For a statement $P(x)$ depending on a point $x\in X$, we say that 
$P(x)$ holds q.e.~$x\in X$ provided the set $\{x\in X:P(x)\text{ holds}\}$ is exceptional (or $m$-polar).\\
 
For the notion of positive continuous additive functional with respect to $\mathbb{X}$, PCAF for short, we adopt the definition from \cite[Definition (1.1), p.~148]{Blumenthal-Getoor} (with $M_t =1_{[0,\zeta)}(t)$, $S=\zeta$).\\
The duality assumption implies  that there is a jointly measurable function 
\[
	p_t(\cdot,\cdot): (0,\infty)\times X\times X\to [0,\infty)
\]
such that 
\[
	\P_t u(x) = \int_X p_t(x,y)u(y)\,dm(y),\quad  t>0,\quad x\in X,\quad u\in\mathscr{B}_b(X).
\]
We assume that $p_t(x,y)>0$ for all $t>0$, $x,y\in X$.\\
For any $\alpha\geq 0$ we designate by $G_\alpha$ the  $\alpha$-order resolvent kernel of $\mathbb{X}$:
\begin{equation}
		G_\alpha (x,y):= \int_0^\infty e^{-\alpha  t} p_t(x,y)\,dt,\quad x,y\in X.
		\label{laplace-transform}
\end{equation}
For $\alp=0$, we set $G:=G_0$.\\
Owing to the  transience assumption we get $0<G_\alpha\not\equiv\infty$. Obviously for any $u\in \mathscr{B}_b(X)$ we have
\[
	R_\alpha u:= \int_X G_\alpha(\cdot,y)u(y)\,dm(y).
\]
We further assume that 
\begin{equation}
		G : X\times X\to (0,\infty]\ \text{is lower semi-continuous  (l.s.c.)}. 
		\label{lsc}
\end{equation}
\begin{rk}
	{\rm We stress that in this section we only need lower semi-continuity of $G(\cdot,y)$ for any fixed $y$. Assumption  (\ref{lsc}) will be crucial in sections 3 and 4. 
		}
\end{rk}
\begin{defi}\label{defi:GKato}
{\rm Let $\mu$ be a positive Radon measure on $X$. 
We say that $\mu$ is
\begin{enumerate}
\item[\rm(1)] a $G$-bounded (Green-bounded) measure if 
\[
G^\mu1:=\int_X G(\cdot,y)\,d\mu(y) \in \mathscr{B}_b(X).
\]
\item[\rm(2)] a $G$-Kato (Green-Kato) measure if   $G^\mu1\in C_0(X)$.
\end{enumerate}
}
\end{defi}
Let $\mathscr{K}(X)$ be the set of $G$-Kato measures.\\
For any  $G$-bounded measure $\mu$ and $u\in\mathscr{B}_b(X)$ we set
\[
	G^\mu u(x) := \int_X G(x,y)u(y)\,d\mu(y),\quad x\in X.
\]
\begin{prop} We have
\begin{enumerate}
\item[\rm(1)] Assume that  $\mu$ is $G$-bounded. Then $G^\mu$ maps continuously $\mathscr{B}_b(X)$ into $\mathscr{B}_b(X)$. If moreover $G^\mu 1$ is continuous then $G^\mu$ maps continuously $\mathscr{B}_b(X)$ into  $C_b(X)$.
\item[\rm(2)] Assume that  $\mu\in  \mathscr{K}(X)$. Then $G^\mu$  maps continuously  $\mathscr{B}_b(X)$ into $C_0(X)$.
\end{enumerate}
\label{PotentialOperator}
\end{prop}	
\begin{proof}
In both cases it is obvious that  $G^\mu:\mathscr{B}_b(X)\to \mathscr{B}_b(X)$ is continuous.\\
Assume that $\mu$ is $G$-bounded and $G^\mu 1$ is continuous. It suffices to prove $G^\mu(\mathscr{B}^+(X))\subset C(X)$. Let $u\in\mathscr{B}^+(X) $, by l.s.c. property for $G$, using Fatou lemma  we get that $G^\mu u$ is l.s.c. On the other hand, by the fact that $G^\mu 1$ is continuous, the first part of the proof and the formula $G^\mu u = \|u\|_{\infty} G^\mu 1 - G^\mu(\|u\|_{\infty} - u)$ we conclude that $G^\mu u$ is upper semi-continuous. Thus $G^\mu u$ is continuous. If furthermore $\mu\in\mathscr{K}(X)$ then the inequality $0\leq  G^\mu u\leq \|u\|_\infty G^{\mu}1$ yields $G^\mu u\in C_0(X)$. 
\end{proof}	
\begin{defi}[{{{See \cite[Th\'eor\`eme VI. 1, p.~524]{Revuz}}}}]\label{defi:smooth}
{\rm Let $\nu$ be a $\sigma$-finite Borel measure. We say that $\nu$ is \emph{smooth (in the strict sense)} if:
\begin{enumerate}
\item[\rm(1)] $\nu$ does not charge any semi-polar sets.
\item[\rm(2)] There exists an increasing sequence of Borel sets 
$\{E_n\}$ such that $X=\bigcup_{n=1}^{\infty}E_n$ and 
\begin{enumerate}
\item[\rm(a)] $\nu(E_n)<\infty$ for each $n\in\mathbb{N}$.
\item[\rm(b)] $\int_{E_n}G_1(x,y)d\nu(y)$ is bounded with respect to $x$.
\item[\rm(c)] Let $T_n:=\sigma_{X\setminus E_n}$ and $T:=\lim_{n\to\infty}T_n$. Then, it should hold 
\begin{align*}
\mathbb{P}_x(T\geq\zeta)=1\quad\text{ for any }\quad x\in X.
\end{align*}
\end{enumerate}
\end{enumerate}
}
\end{defi}

\begin{prop}\label{prop:Smooth}
\begin{enumerate}	
\item[\rm(1)] Any $G$-Kato measure is smooth in the sense of Definition~\ref{defi:smooth}.
\item[\rm(2)] Any $G$-bounded measure charging no semi-polar sets is smooth in the sense of Definition~\ref{defi:smooth}. 
\end{enumerate}
\end{prop}
\begin{proof}
{\em Step 1:} Let $\nu$ be a $G$-Kato measure. To prove the first property from Definition~\ref{defi:smooth}, we may assume $\nu(X)<\infty$. Otherwise, 
 taking an increasing sequence $\{A_n\}$ satisfying $\nu(A_n)<\infty$ 
 for each $n\in\mathbb{N}$ and $X=\bigcup_{n=1}^{\infty}A_n$; then 
 $1_{A_n}\nu$ charges no semi-polar set for each $n$, hence $\nu$ charges no semi-polar set.\\
 For $\alpha\geq 0$, let $R_{\alpha}\nu(x):=\int_XG_{\alpha}(x,y)d\nu(y)$. 
 Set $u:= G^\nu 1$ and for all $n\in\N$,  $u_n:= nR_n u$. Then $0\leq u_n\leq u$. Indeed, from the resolvent identity we get
 \begin{align*}
u_n =  nR_n u(x) &= n\int_X G_n(x,y) G^\nu 1(y)\,dm(y) = n\int_X \big( \int_X G_n(x,y) G(y,z)\,dm(y) \big) d\nu(z)\\
 & \int_X ( G(x,z) - G_n(x,z) d\nu(z)\leq 	\int_X  G(x,z) d\nu(z) =u.
 \end{align*}
 Let us show that $(u_n)_n$ is increasing and $\lim_{n\to\infty} u_n =u$, pointwise. Let $m\geq n$, using the resolvent identity once again and the inequality $0\leq u_n\leq u$ we obtain
 \begin{align*}
 u_n &= nR_n u = nR_m u + n(m-n) R_m R_n u = nR_m u + (m-n) R_m u_n\\
 & \leq n R_m u + (m-n) R_m u = m R_m u = u_m.
 \end{align*}
 Now using the relationship between the resolvent and the semigroup we conclude that $u_n$ increases to $u$ point-wisely.\\
 Since the measure $m$ is $G_n$-bounded and $R_n1$ is continuous  and $u$ is in $C_0(X)$ we get, by Proposition \ref{PotentialOperator}, $(u_n)\subset C_b(X)$ and hence $(u_n)_n\subset C_0(X)$, by $0\leq u_n\leq u$. Using Dini's lemma we conclude that  $u_n\uparrow u$ uniformly. Since $R_{\alpha}\nu$ is bounded and $\nu(X)<\infty$, 
 in the same way of the first paragraph in \cite[the proof of Theorem~6.4.5]{Oshima}, there exists a positive continuous additive functional ${\bf A}^{\nu}$ admitting exceptional set associated to $\nu$ such that 
 \begin{align*}
 	R_1\nu(x)=\mathbb{E}_x\left[\int_0^{\infty}e^{-t}d {\bf A}^{\nu}_t \right] \quad\text{ q.e.}\quad x\in X.
 \end{align*}
This correspondence remains valid 
 between general smooth measures and PCAF also 
 by \cite[Th\'eor\`eme VI. 1, p.~524]{Revuz}.
 Thanks to the resolvent equation: 
 $R_{\alpha}\nu(x)-R_{\beta}\nu(x)+(\alpha-\beta)R_{\alpha}R_{\beta}\nu(x)=0$, we can deduce that for any $\alpha>0$
 \begin{align*}
 	R_{\alpha}\nu(x)=\mathbb{E}_x\left[\int_0^{\infty}e^{-\alpha t}d {\bf A}^{\nu}_t \right]\quad\text{ q.e.}\quad x\in X.
 \end{align*}
 Similarly, for any bounded measurable function $f$, we have 
 \begin{align*}
 	R_{\alpha}(f\nu)(x)=\mathbb{E}_x\left[\int_0^{\infty}e^{-\alpha t}f(X_t)d {\bf A}^{\nu}_t \right] \quad\text{ q.e.}\quad x\in X.
 \end{align*}
 If $B$ is a nearly Borel semi-polar set, then by \cite[Proposition 3.4, p.~80]{Blumenthal-Getoor} $\mathbb{X}_t\in B$ for only countably many values of $t$, $\mathbb{P}_x$ almost surely for all $x\in X$. Thus using continuity of ${\bf A}^\nu$ we obtain for any $\alpha>0$, 
 \begin{align*}
 	R_{\alpha}(1_B\nu)(x)=\mathbb{E}_x\left[\int_0^{\infty}e^{-\alpha t}1_B(X_t)d {\bf A}^{\nu}_t \right]=0 \quad\text{ q.e.}\quad x\in X.
 \end{align*}
 Hence 
 \begin{align*}
 	0=\lim_{\alpha\to\infty}\alpha \int_XfR_{\alpha}1_B\nu\, dm=
 	\lim_{\alpha\to\infty}\langle 1_B\nu, \alpha R_{\alpha}^*f\rangle=
 	\langle 1_B\nu,f\rangle
 \end{align*}
 for any $f\in C_c(X)_+$, where $R_\alp^*$ is the dual resolvent. Thus $\nu(B)=0$.\\ 
 
{\em Step 2}: Let $\nu$ be either a $G$-Kato measure or  a $G$-bounded measure charging no semi-polar sets. It suffices to prove (2c). The other conditions (2a) and (2b) are easily confirmed. Next we prove (2c). Let $\{E_n\}$ be an increasing sequence of relatively compact open subsets of $X$ satisfying $\overline{E}_n\subset E_{n+1}$ for all $n\in\mathbb{N}$ and $X=\bigcup_{n=1}^{\infty}E_n$.  Let $\varphi$ be a strictly positive 
bounded measurable function with $\int_X\varphi\, dm=1$. We now define an outer capacity $\Gamma$ on $X$ as follows: 
\begin{align*}
\Gamma(O)&:=\mathbb{E}_{\varphi m}[e^{-\sigma_{O}}]\quad \text{ for an open subset }O\text{ of }X,\\
\Gamma(A)&:=\inf\{\Gamma(O): A\subset O, \;O\text{ is open}\}\quad\text{ for a subset }A\text{ of }X.
\end{align*}
It is clear that $\Gamma(\emptyset)=0$, $\Gamma(X)=1$ and $0\leq\Gamma(A)\leq1$ for any $A\subset X$. 
Moreover, $\Gamma$ is a Choquet capacity: 
\begin{enumerate}
\item[\rm (i)] For any subsets $A_1,A_2$ of $X$ with 
$A_1\subset A_2$, we have  
$\Gamma(A_1)\leq\Gamma(A_2)$. 
\item[\rm (ii)] Let $\{A_n\}$ be an increasing sequence of subsets of $X$ and $A:=\bigcup_{n=1}^{\infty}A_n$. Then $\Gamma(A)=\lim_{n\to\infty}\Gamma(A_n)$. 
\item[\rm (iii)] Let $\{K_n\}$ be an decreasing sequence of compact subsets of $X$ and $K:=\bigcap_{n=1}^{\infty}K_n$. Then $\Gamma(K)=\lim_{n\to\infty}\Gamma(K_n)$.
\end{enumerate}
To prove (i)--(iii), we confirm only the following strong subadditivity of $\Gamma$: 
\begin{enumerate}
\item[\rm (iv)] For any open subsets $O_1,O_2$ of $X$, 
\begin{align*}
\Gamma(O_1\cup O_2)+\Gamma(O_1\cap O_2)\leq\Gamma(O_1)+\Gamma(O_2).
\end{align*}
\end{enumerate}
Since $\sigma_{O_1\cup O_2}=\sigma_{O_1}\land\sigma_{O_2}$ and 
$\sigma_{O_1\cap O_2}\geq\sigma_{O_1}\lor\sigma_{O_2}$, (iv) is easily confirmed in view of the elementary equality: 
$e^{-a\land b}+e^{-a\lor b}=e^{-a}+e^{-b}$ for $a,b\geq0$.
Then we can confirm (i)--(iii) in the same way of the proof of \cite[Theorem~1.2.10]{CF}.\\ 
Moreover, our capacity $\Gamma$ is tight in the sense that 
there exists an increasing sequence $\{K_n\}$ 
of compact sets such that 
\begin{align}
\lim_{n\to\infty}\Gamma(X\setminus K_n)=\lim_{n\to\infty}\mathbb{E}_{\varphi m}[e^{-\sigma_{X\setminus  K_n}}]=
0.\label{eq:tight}
\end{align}
Indeed, since the probability measure $\mathbb{P}_{\varphi m}$ on $\mathscr{M}$ is complete and the sample paths of $(\mathbb{X}_t)$ is 
c\`adl\`ag, there exists an increasing sequence $\{K_n\}$ 
of compact sets such that $\mathbb{P}_{\varphi m}(\lim_{n\to\infty}\sigma_{X\setminus K_n}<\zeta)= \mathbb{P}_{\varphi m}(\lim_{n\to\infty}\sigma_{X\setminus K_n}<\infty)=0$ by \cite[Chapter IV, alternative proof of Theorem~1.15]{Ma-Rockner}. 
Then we have \eqref{eq:tight}.\\
In view of the tightness of 
$\Gamma$, any closed set $F$ is quasi-compact in the sense of 
\cite[Definition~2.1 and Lemma~2.2]{Fuglede}. 
Combining this with \cite[Theorem~2.10]{Fuglede}, $\Gamma$ is continuous below for decreasing sequence of closed sets 
in the sense that for any decreasing sequence of closed subsets $\{F_n\}$ of $X$, we have 
\begin{align}
\lim_{n\to\infty}\Gamma(F_n)=\Gamma\left(\bigcap_{n=1}^{\infty}F_n \right).\label{eq:decreasing}
\end{align}
Next we show that for any closed subset $F$ of $X$
\begin{align}
\Gamma(F)\geq\mathbb{E}_{\varphi m}[e^{-\sigma_{F}}].\label{eq:closedCompa}
\end{align}
Let $O_n:=\{x\in X: d(x,F)<1/n\}$. Then $A_n:=\{x\in X: d(x,F)\leq1/n\}$ forms a decreasing sequence of closed sets containing $F
=\bigcap_{n=1}^{\infty}A_n$. Here $d$ is the distance such that 
any open $B_r(x)=\{y\in X:d(x,y)<r\}$ is relatively compact for $r>0$. Applying the continuity below for the decreasing sequence of closed sets, we have 
\begin{align*}
\Gamma(F)&=\lim_{n\to\infty}\Gamma(A_n)
=\lim_{n\to\infty}\Gamma(O_n)\\
&=\lim_{n\to\infty}\mathbb{E}_{\varphi m}[e^{-\sigma_{O_n}}]
=\mathbb{E}_{\varphi m}[e^{-\lim_{n\to\infty}\sigma_{O_n}}]\geq \mathbb{E}_{\varphi m}[e^{-\sigma_F}],
\end{align*}
where we use $\sigma_F\geq \lim_{n\to\infty}\sigma_{O_n}$. 
Apply \eqref{eq:decreasing} for $F_n:=X\setminus E_n$ so that 
$\lim_{n\to\infty}\Gamma(F_n)=\Gamma(\emptyset)=0$. Then by \eqref{eq:closedCompa}, 
we have $\lim_{n\to\infty}\mathbb{E}_{\varphi m}[e^{-\sigma_{X\setminus E_n}}]=0$, hence $\mathbb{P}_x(
T=\infty)=1$ holds $m$-a.e.~$x\in X$ for $T:=\lim_{n\to\infty}\sigma_{X\setminus E_n}$. Finally, we prove 
$\mathbb{P}_x(T=\infty)=1$ for all $x\in X$. For any $x\in X$, we have 
\begin{align*}
0&=\int_X G_2(x,y)\mathbb{E}_y[e^{-T}]dm(y)\\
&=R_2(\mathbb{E}_{\cdot}[e^{-T}])(x)\\&=\int_0^{\infty}
e^{-2t}{\bf P}_t(\mathbb{E}_{\cdot}[e^{-T}])(x) dt\\
&=\int_0^{\infty}e^{-t}\mathbb{E}_x[e^{-(t+T\circ \theta_t)}]dt,
\end{align*}
which implies $\mathbb{E}_x[e^{-(t+T\circ \theta_t)}]=0$ for all  $t>0$, because $t\mapsto t+T\circ \theta_t$ is right continuous at $t\in(0,\infty)$. 
In particular, $\mathbb{E}_x[e^{-T}:t<T]=0$ for all~$t>0$, hence 
$\mathbb{E}_x[e^{-T}:0<T]=0$.  
Then 
\begin{align*}
\mathbb{E}_x[e^{-T}]&=\mathbb{E}_x[e^{-T}:T>0]+\mathbb{E}_x[e^{-T}:T=0]\\
&=\mathbb{P}_x(T=0)\\
&=\mathbb{P}_x(\sigma_{X\setminus E_n}=0\text{ for all }n\in\mathbb{N})\quad (\because\sigma_{X\setminus E_n}\leq T)\\
&=\downarrow\lim_{n\to\infty}\mathbb{P}_x(\sigma_{X\setminus E_n}=0)
=\downarrow\lim_{n\to\infty}1_{(X\setminus E_n)^r}(x)\\
&\leq \downarrow\lim_{n\to\infty}1_{X\setminus E_n}(x)=1_{\bigcap_{n=1}^{\infty}(X\setminus E_n)}(x)=0.
\end{align*}
This implies $\mathbb{P}_x(T=\infty)=1$. 
Here we use the Blumenthal's zero-one law and 
$(X\setminus E_n)^r\subset(X\setminus E_n)$. 
\end{proof}
\begin{rk}
{\rm  
\begin{enumerate}	
\item[\rm(1)]In the framework of (non-symmetric) regular Dirichlet forms, 
any finite $G$-bounded measure is of finite energy integral, in particular, it is smooth in the sense of \cite[(4.5)]{Ma-Rockner}.  Hence the assertion of Proposition~\ref{prop:Smooth} is consistent with this, because any $G$-bounded measure does not charge $m$-polar sets and any semi-polar sets are $m$-polar in the framework of semi-Dirichlet forms (see \cite{Silverstein}). 
\item[\rm(2)] Proposition \ref{prop:Smooth} still holds true if we replace $G$-Kato, respectively $G$-bounded by $G_\alp$-Kato, respectively $G_\alp$-bounded for some $\alp>0$. It also holds true if we replace $\nu$ is $G$-Kato by $\varphi\cdot\nu$ is $G_\alp$-Kato for some $\alp>0$ and some $\varphi\in C_0(X)$ with $\varphi>0$.
\item[\rm(3)] Assertion (2) from  Proposition \ref{prop:Smooth} can be derived from \cite[Theorem 3.13, p.~167 and Theorem 3.5, p.~281]{Blumenthal-Getoor} in conjunction with \cite[Th\'eor\`eme VI. 1, p.~524]{Revuz}. However, we give here a new proof.
\end{enumerate}
}
\end{rk} 
Let $\mu$ be a Borel positive smooth measure and let ${\bf A}$ be the (unique up to equivalence class) PCAF whose Revuz measure is $\mu$ with respect to the process $\mathbb{X}$. The (Revuz) correspondence between $\mu$ and $\bf{A}$ is expressed through the formula (see \cite[Proposition, p.~507]{Revuz}): for $u\in\mathscr{B}^+(X)$ it holds
\begin{equation}
\int_X u\,d\mu = \lim_{t\downarrow 0} \frac{1}{t}\mathbb{E}_m\left[\int_0^t u(\mathbb{X}_s)\,d{\bf A}_s\right]
=\lim_{\alp\to\infty} \alp\int_X U_\alp^{\bf A} u(x)\,dm(x),\quad u\in\mathscr{B}^+(X),
\label{Revuz}
\end{equation}	
where
\[
U_\alp^{\bf A} u(x) = \mathbb{E}_x\left[ \int_0^\infty e^{-\alp t} u(\mathbb{X}_t)\,d{\bf A}_t \right],\quad x\in X,
\]
is the $\alp$-potential of $u$ with respect to the PCAF $\bf{A}$. Furthermore, by \cite[{Th\'eor\`eme, p.~522}]{Revuz} we have
\begin{equation}
U_\alp^{\bf A} u(x) = G_\alp^\mu u (x):= \int_X G_\alp (x,y) u(y)\,d\mu(y),\quad u\in\mathscr{B}^+(X),\quad \alp\geq 0.
\label{Apotential-Mpotential}
\end{equation}
If moreover $\mu$ is  $G$-bounded, then both  identities (\ref{Revuz}), (\ref{Apotential-Mpotential}) extend to every $u\in\mathscr{B}_b(X)$.\\
Let us stress that, by Proposition \ref{prop:Smooth} any $G$-Kato measure or any $G$-bounded measure charging no semi-polar set has a unique PCAF in the sense precised above.\\ 
Let  $\tau_t, t\geq 0$ be the right-continuous inverse of ${\bf A}_t$:
\[
\tau_t = \inf\left\{s\geq0: {\bf A}_s>t\right\}.
\]
Let us consider the time-changed process
\[
\check{\mathbb{X}}:= \left( \Om,\mathscr{M},\mathscr{M}_{\tau_t}, (\mathbb{P}_x)_{x\in X_{\partial}}, (\mathbb{X}_{\tau_t})_{t\geq 0}, (\theta_{\tau_t})_{t\geq 0}\right),
\]
with state space $(X,\mathscr{B}(X))$ augmented by $\partial$.\\
It is known that the process $\check{\mathbb{X}}$ considered on $X$ is a strong Markov process. At this stage we refer to 
 \cite[pp.~212--213]{Blumenthal-Getoor}. However, in general it is not a standard process because of the lack of normality.\\
The semigroup and the resolvent of $\check{\mathbb{X}}$ are given respectively by: For any $t> 0, \alpha\geq 0$
\[
\check{\bf P}_t u(x) := \mathbb{E}_x\left[u(\mathbb{X}_{\tau_t}) \right],\  
\check{R}_\alpha u(x) := \mathbb{E}_x\left[\int_0^\infty e^{-\alpha t}u(\mathbb{X}_{\tau_t})\,dt \right], \quad 
 u\in\mathscr{B}_b(X),\quad x\in X.
\] 
\begin{rk}
{\rm
We mention that Revuz correspondence between certain class of measures and PCAF has been  established in a more general context by Fitzsimmons--Getoor in \cite{Fitzsimmons-Getoor96} and Beznea--Boboc in \cite{Beznea}. Also time-change transformation can be made for more general processes as in  \cite{Sharpe}.  
}
\end{rk}
The content of the following proposition is known; for the convenience of the reader we include it together with its proof.
\begin{prop}
\begin{enumerate}
\item[\rm(1)] The family $(\check{R}_\alp)_{\alp\geq 0}$ satisfies the resolvent equation: $\check{R}_\alp - \check{R}_\beta = (\beta- \alp)\check{R}_\alp\check{R}_\beta, \alp,\beta\geq 0$. Moreover $\alp\|\check{R}_\alp\|\leq 1, \alp\geq 0$.
\item[\rm(2)] For all $\alp\geq 0, \ran(\check{R}_\alp)$ and $\ker(\check{R}_\alp)$ are $\alp$-independent.
\item[\rm(3)]   For all $\alp\geq 0, \check{R}_\alp$ is injective on $\ran(\check{R}_\beta)$ for all $\beta\geq 0$.
\end{enumerate}
\label{First-Prop-Res}
\end{prop} 
\begin{proof}
(1): Let $u\in\mathscr{B}_b(X)$ and  $0\leq\beta<\alp$. Then, by Markov property we get
\begin{align*}
\check{R}_\alp \check{R}_\beta u(x) &= \int_0^{\infty} e^{-\alp t} \check{\bf P}_t(\check{R}_\beta u)(x)\,dt
= \int_0^{\infty} e^{-\alp t} \check{\bf P}_t\left(\int_0^\infty e^{-\beta s} \check{{\bf P}}_s u(x)\,ds\right)\,dt\\
& = \int_0^{\infty} \int_0^{\infty} e^{-\alp (t + s)} \check{\bf P}_t \check{{\bf P}}_s u(x)\,ds\,dt
= \int_0^{\infty} \int_0^{\infty} e^{-\alp t}  e^{-\beta s}\check{\bf P}_{t+s} u(x)\,ds\,dt\\
& = \int_0^\infty e^{(\beta - \alp)t}\left(\int_0^\infty  e^{-\beta s}\check{\bf P}_{s} u(x)\,ds\right)\,dt 
- \int_0^\infty  e^{(\beta - \alp) t}\left( \int_0^t  e^{-\beta s}\check{\bf P}_{s} u(x)\,ds\right)\,dt\\
& = -\frac{1}{\beta - \alp} \check{R}_\beta u(x) -  
\int_0^\infty e^{-\beta s}\check{\bf P}_{s} u(x) \left(\int_s^\infty  e^{(\beta - \alp)t}\,dt\right)\,ds\\
& = -\frac{1}{\beta - \alp} \check{R}_\beta u(x)  + \frac{1}{\beta - \alp} \check{R}_\alp u(x), 
\end{align*}
which yields the resolvent equation.\\
(2): Follows directly from the resolvent equation.\\
(3): Let $u=\check{R}_\beta v\in\ker(\check{R}_{\alp_0})$. By assertion (2) we get that $u\in\ker(\check{R}_{\alp})$ for all $\alp>0$.  From the resolvent identity we get 
\[
0=\alp\check{R}_\alp \check{R}_\beta v = \alp\check{R}_\beta\check{R}_\alp v = \frac{\alp}{\beta - \alp}\check{R}_\alp v - \frac{\alp}{\beta - \alp}\check{R}_\beta v.
\]
Thus letting $\alp\to\infty$ we get $\check{R}_\beta v =0$ and hence $\check{R}_{\alp_0}$ is injective on $\ran(\check{R}_\beta)$ for all $\beta\geq 0$.
\end{proof}
\begin{prop}
For every  $u\in\mathscr{B}_b(X)$ we have
\begin{equation}
	\check{R}_\alpha u(x) = \mathbb{E}_x\left[\int_0^\infty e^{-\alpha {\bf A}_t}u(\mathbb{X}_t)\,d{\bf A}_t \right],\quad \alp>0.
\end{equation}
Assume in addition that $\mu$ is $G$-bounded and charges no semi-polar sets. Then 
\begin{equation}
	\check{R}_\alpha = (1 + \alpha G^\mu)^{-1} G^\mu,\quad \alpha\geq0.	
\label{Resolvent-Formula}	
\end{equation}
\label{ResolventFormula}
\end{prop}
\begin{proof}
The first formula follows from the substitution $t\rightarrow {\bf A}_t$ (see \cite[p.~207]{Blumenthal-Getoor}).\\	
Assume that $\mu$ is $G$-bounded. Let $u\in\mathscr{B}_b(X)$. From the first formula together with the identity   (\ref{Apotential-Mpotential}) we infer  that $\check{R}_0 u = G^\mu u$; and hence formula (\ref{Resolvent-Formula}) is true for $\alp=0$. This leads by the resolvent equation to $( 1+ \alp G^\mu) \check{R}_\alp= G^\mu $. Let us show that the operator  $1 + \alp G^\mu$ is injective for all $\alp$. For, let $u\in\mathscr{B}_b(X)$ be such that  $u + \alp G^\mu u = 0$. Then $0 = \check{R}_\alp ( 1+ \alp G^\mu)u = G^\mu u$. Owing to $u + \alp G^\mu u = 0$ we obtain $u=0$. Thus the operator $1+ \alp G^\mu$ is injective and then  $\check{R}_\alp =  ( 1+ \alp G^\mu)^{-1} G^\mu$.
\end{proof}
Applying Proposition \ref{PotentialOperator} in conjunction with Proposition \ref{ResolventFormula} we get:
\begin{coro}
\begin{enumerate}
\item[\rm(1)] Assume that $\mu$ is $G$-bounded, charges no semi-polar sets and $G^\mu 1$ is continuous. Then  $\check{R}_\alpha$ maps continuously $\mathscr{B}_b(X)$ into $C_b(X)$.	
\item[\rm(2)] Assume that $\mu\in\mathscr{K}(X)$. Then $\check{R}_\alpha$ maps continuously $\mathscr{B}_b(X)$ into $C_0(X)$.
\end{enumerate}
\label{Feller-Resolvent}
\end{coro}
\begin{rk}
{\rm 
We mention that compared to results from \cite[Theorem 7.1]{Kuwae-17}, Corollary  \ref{Feller-Resolvent} needs no extra assumption on the starting process $\mathbb{X}$ such as symmetry or the doubly Feller property for its resolvent. However, the $G$-Kato assumption on the measure $\mu$ is stronger than the local Kato  assumption used in \cite[Theorem 7.1]{Kuwae-17}.
}
\end{rk}
The following corollary has an independent interest. For instance its content is extensively used in \cite{Hansen-Global} to get, among other results, comparison for perturbed Green functions.\\
For the definition of the \lq complete maximum principle\rq\, (used below) we refer the reader to \cite[p.~76]{Bliedtner-Hansen}.
\begin{coro}
For every $G$-bounded measure $\mu$ charging no semi-polar sets,  the kernel $G^\mu$ satisfies the complete maximum principle.
\end{coro}
\begin{proof}
On the one hand  $\check{R}_\alpha\to G^\mu$ as $\alpha\to 0$, in the operator norm. On the other one, $\alpha \check{R}_\alpha$ is sub-Markovian, which is equivalent, by \cite[Theorem 7.7, p.~79]{Bliedtner-Hansen} to the fact that $G^\mu$ satisfies the complete maximum principle.
\end{proof}
We already mentioned that $\check{\mathbb{X}}$ considered on $X$ needs not to be a standard process due to a lack of normality. Even more, $\check{{\bf P}}_t$ is degenerate and $\check{R}_\alp$ is not injective. In fact, for all $u\in\mathscr{B}_b(X)$ such that $u$ vanishes on the complement of the support of $\mu$ it holds $\check{{\bf P}}_t u = \check{R}_\alpha u = 0$ for all $t,\alp\geq 0$.\\ 
We turn our attention now to give  a necessary and sufficient condition for $\check{\mathbb{X}}$ to obey the normality property and hence to be a standard process on $X$. To this end we introduce the fine support of a Borel smooth measure $\mu$; it is the smallest nearly Borel measurable finely closed set $B\subset X$ such that $\mu(B^c)=0$. Whereas the fine support  of $\bf{A}$ is the smallest nearly Borel measurable finely closed set $B\subset X$ such that $U_0^{\bf A}(1_{B^c})=0$, where $U_0^{\bf A}$ is the potential of $\bf{A}$. Let us quote that by our assumptions the measure $m$ is a reference measure, as defined in \cite[Definition (1.1), p.196]{Blumenthal-Getoor}, and hence the fine support is well defined. 
\begin{prop}
The fine support of $\mu$ coincides with the fine support of $\bf{A}$.
\label{both-supports}
\end{prop}
\begin{proof}
The proof follows from  formula (\ref{Apotential-Mpotential}) which  yields $U_0^{\bf A}(1_{B^c}) = G^\mu 1_{B^c}$ in conjunction with the fact that $G>0$.
\end{proof}
Set 
\[
\varphi^{\bf A}(x):= \mathbb{E}_x(e^{-\tau_0}),\quad x\in X.
\]
By \cite[Corollary (3.10), p.~215]{Blumenthal-Getoor} together with Proposition \ref{both-supports}, both fine supports coincide with the set
\begin{align*}
F  : &= \left\{ x\in X: \varphi^{\bf A}(x) =1\right\}\\
& = \left\{ x\in X: \mathbb{P}_x\{\tau_0=0\} =1\right\}.
\end{align*}
\begin{prop} The following assertions hold:
\begin{enumerate}	
\item[\rm(1)] The set $F$ is Borel measurable.
\item[\rm(2)] The process $\check{\mathbb{X}}$ with state space $F$ is a Borel right-continuous Markov process.
\end{enumerate}
\label{F-Borel}
\end{prop}
\begin{proof}
(1): Owing to the fact that $m$ is a reference measure, we  deduce from \cite[Proposition 1.4, p.~197]{Blumenthal-Getoor} that $\varphi^{\bf A}$ is Borel measurable and hence $F$ is a Borel set.\\
(2): Follows from (1) and \cite[Theorem 65.9, p.~307]{Sharpe}.
 \end{proof}
\begin{rk}
{\rm
In many occasions later we shall assume that $F$ is closed. We claim that if $F$ is closed then it coincides with the topological support of $\mu$. Indeed, let $F_{\rm top}$ be the topological support of the measure $\mu$. Since the fine topology is finer than the original topology we get $F\subset F_{\rm top}$. Hence if moreover $F$ is closed it follows from the definition of the fine and the topological supports that $F_{\rm top}\subset F$ and then $F_{\rm top}= F$.
}
\end{rk}
Let $\sigma_{F}$ be the hitting time of $ F$, i.e.,
\[
	\sigma_F = \inf\{t>0: \mathbb{X}_t\in F \}.
\]
By \cite[Proposition 3.5, p.~213]{Blumenthal-Getoor} we have	
\[
\mathbb{P}_x\{ \tau_0= \sigma_F\} = 1,\quad x\in X.
\]
Let $P_F$  be the hitting operator
\begin{equation}
	P_F u(x) := \mathbb{E}_x\left[u(X_{\sigma_F})\right],\quad u\in\mathscr{B}_b(X_\partial)\cup\mathscr{B}^+(X_\partial),\quad x\in X_\partial.
	\label{hittin-op1}
\end{equation}
Then 
\begin{equation}
	P_F u(x) =  \mathbb{E}_x\left[u(\mathbb{X}_{\tau_0})\right],\quad  u\in\mathscr{B}_b(X)\cup\mathscr{B}^+(X),\quad x\in X.
	\label{hitting-op2}
\end{equation}
According to \cite[p.~213]{Blumenthal-Getoor} we have
\begin{equation}
	\mathbb{P}_x\left\{\mathbb{X}_{\tau_0} = x  \right\} =1 \Leftrightarrow x\in F,
	\label{normality}
\end{equation}
which leads in conjunction with (\ref{hitting-op2}) to
\begin{equation}
	P_F u = u\quad \text{on}\quad F,\quad u\in\mathscr{B}_b(X)\cup\mathscr{B}^+(X).
	\label{regularisee}	
\end{equation}
\begin{lem}
Let $u\in C_b(X)$. Then 
\begin{equation}
		\lim_{t\downarrow 0} \check{\bf P}_t u(x) =  \lim_{\alpha\to\infty} \alpha \check{R}_\alpha u(x)  
		= \mathbb{E}_x\left[u(\mathbb{X}_{\tau_0})\right] = P_F u(x),\quad x\in X.
\end{equation}
\label{sg-res-limit}
\end{lem}
\begin{proof}
Let $u\in C_b(X)$. The right continuity of $\mathbb{X}_{\tau_t}$ in conjunction with dominated convergence theorem  leads to $\lim_{t\downarrow 0} \check{\bf P}_t u(x) = \mathbb{E}_x\left[u(\mathbb{X}_{\tau_0})\right],\; x\in X$. Besides
\[
	\alpha \check{R}_\alpha u(x) = \int_0^\infty e^{-t}\check{{\bf P}}_{\frac{t}{\alp}}u(x)\,dt\to \mathbb{E}_x\left[u(\mathbb{X}_{\tau_0})\right],\quad x\in X\quad \text{as}\quad \alp\to\infty,
\]
by the first part of the proof and dominated convergence theorem.
\end{proof}
Lemma \ref{sg-res-limit} in conjunction with the identity (\ref{regularisee}) leads to
\begin{equation}
	\lim_{t\downarrow 0} \check{\bf P}_t u(x) =  \lim_{\alpha\to\infty} \alpha \check{R}_\alpha u(x)  
	= u(x), \quad x\in F,\quad u\in C_b(X).
	\label{conv-fine-supp}	
\end{equation}
We now turn our attention to give necessary and sufficient conditions for $\check{\mathbb{X}}$ to be a standard process on $X$.
\begin{theo} The following assertions are equivalent.
\begin{enumerate}
\item[\rm(1)] The process $\check{\mathbb{X}}$ is a standard process on $X$.
\item[\rm(2)] $F=X$.
\item[\rm(3)] $P_F u = u$ for all $u\in \mathscr{B}_b(X)$.
\item[\rm(4)] $P_F$ is injective on $\mathscr{B}_b(X)$.
\end{enumerate}
\label{right}
\end{theo}
\begin{proof}
We recall that $\check{\mathbb{X}}$ satisfies all the conditions  of the standard process but not the normality.\\
(1)$\Leftrightarrow $(2): 
Assume that $F=X$. From the normality of $\mathbb{X}$ we get $\mathbb{P}^\partial\{X_{\tau_t} = \partial \} =1$,  for all $t\geq 0$.\\
Let $x\in X$, by assumptions and from the definition of the fine support we get $\mathbb{P}_x\{ \tau_0 = 0 \} =1$. Thus, using normality of $\mathbb{X}$ once again, we get
\[
\mathbb{P}_x \{ X_{\tau_0} = x \} = \mathbb{P}_x \{ X_0= x \} = 1,\quad x\in X,
\]
and then $\check{\mathbb{X}}$ is normal.\\
Conversely: Assume that $\mathbb{P}_x \{ X_{\tau_0} = x \}=1$ for all $x\in X_\partial$. From the property $\mathbb{P}_x \{ X_{\tau_t}\in F_\partial\ \text{for all}\ t\geq 0 \} = 1$ for all  $x\in X_\partial$ (see \cite[Theorem 65.9, p.~307]{Sharpe}), we get that  $\mathbb{P}_x \{ X_{\tau_0}\in F_\partial \} = 1,\, x\in X_\partial$. In particular for $x\in X$ we get  $\mathbb{P}_x \{ X_{\tau_0}\in F \} = 1$ and hence $x\in F$. Thus $X\subset F$ and hence $X=F$.\\
(2)$\Rightarrow$(3) is obvious.\\
(3)$\Rightarrow$(2). Since by Proposition \ref{F-Borel} the set $F$ is Borel measurable we can take $u=1_F$ under (3), then we get $F=X$.\\
(3)$\Rightarrow$(4) is obvious.\\
(4)$\Rightarrow$(3): By the resolvent identity (or Markov property for $\check{\bf P}_t$) we get  $P_F u = P_F(P_F u)$ for all $u\in C_b(X)$ and hence for all $u\in\mathscr{B}_b(X)$ by monotone class argument. Thus $P_F(u-P_F u)=0$ and by assumption we get $u-P_Fu =0$.
\end{proof}
\begin{rk}
{\rm
A sufficient condition ensuring one of the equivalent conditions of Theorem \ref{right} to be true is that ${\bf A}_t$ is strictly increasing.
}
\end{rk}
From now on we assume that $\mu\in\mathscr{K}(X)$. 
\begin{lem}
\label{Feller-1}	
Assume that $F$ is closed. Then
\begin{enumerate}
\item[\rm(1)] $P_F(C_0(X))\subset C_0(X)$.
\item[\rm(2)] For all $u\in C_0(X)$ we have $\lim_{\alpha\to\infty}\alpha \check{R}_\alpha u = P_Fu$, uniformly on $X$. 
\end{enumerate}
\end{lem}
\begin{proof}
To prove the first assertion we will use \cite[Th\'eor\`eme IV, p.~425]{Lion}. However, as $X$ is not assumed to be compact we have to go over to the one point compactification of $X$.\\
We mention that the semigroup ${\bf P}_t$ is Markovian on the space $\mathscr{B}_b(X_\partial)$ and hence the semigroup $\check{{\bf P}}_t, t>0$ defines a Markovian semigroup on the space $\mathscr{B}_b(X_\partial)$: for $u\in\mathscr{B}_b(X_\partial),\, x\in X_\partial$ we have
\[
\check{{\bf P}}_t u(x) = \mathbb{E}_x\left[u(\mathbb{X}_{\tau_t})\right], \quad
\check{{\bf P}}_t 1_{X_\partial}(x) = \mathbb{E}_x\left[1_{X_\partial}(\mathbb{X}_{\tau_t})\right](x) = 1,\quad {x\in X_{\partial}}.
\]
Hence for any $\alpha$ the resolvents $\check{R}_\alpha$ are well defined on $\mathscr{B}_b(X_\partial)$ and still verify the resolvent identity. Moreover for all $\alpha>0$,  $\check{R}_\alpha$ is Markovian:
\[
\alp\check{R}_\alpha 1_{X_\partial} = 1\quad {\text{ on }\quad X_{\partial}}. 
\]
Let us show that $\check{R}_\alpha$ maps $C(X_\partial)$ into $C(X_\partial)$. Let $u\in C(X_\partial)$. Then
\begin{align*}
\check{R}_\alpha u(x)& = \check{R}_\alpha (u|_X)(x) + \check{R}_\alpha (u|_{\{\partial\}})(x) 
= \check{R}_\alpha (u|_X)(x) + u(\partial)\check{R}_\alpha (1_{\{\partial\}})(x). 
\end{align*}
Let us compute $\check{R}_\alpha (1_{\{\partial\}})(x), x\in X_\partial$. By Markov property, for any $x\in X_\partial$, we have 
\begin{align*}
\check{R}_\alpha 1_{X_\partial}(x)	= \frac{1}{\alp} = \check{R}_\alpha (1_X)(x) + \check{R}_\alpha (1_{\{\partial\}})(x). 
\end{align*}
For $x\in X$ we have $\check{R}_\alpha (1_X)(x)  = \check{R}_\alpha 1(x)$ (the resolvent on $C_b(X)$ evaluated at the constant $1$) and for $x=\partial$ we have
\begin{align*}
\check{R}_\alpha (1_X)(\partial) = \mathbb{E}_\partial\left[ \int_0^\infty {e^{-\alpha t}}1_X(\mathbb{X}_{\tau_t})\,dt\right].
\end{align*}
Since $\mathbb{E}_\partial \{ \mathbb{X}_t\in X\} = 0$ for all $t\geq 0$ we get $\check{R}_\alpha (1_X)(\partial) =0$. Thus 
\begin{align*}
\check{R}_\alpha (1_{\{\partial\}})(x)	= \frac{1}{\alp} - \check{R}_\alpha (1_X)(x) 
=
\begin{cases} 
\quad \frac{1}{\alp} - \check{R}_\alpha 1(x),\quad x\in X, \\
\quad \frac{1}{\alp}, \quad x =\partial. 	
\end{cases}
\end{align*}
Similar computations lead to 
\begin{align*}
\check{R}_\alpha u(x)
=
\begin{cases} 
\quad \check{R}_\alpha (u|_X))(x)  +   \frac{u(\partial)}{\alp} - u(\partial)\check{R}_\alpha 1(x),\quad x\in X, \\
		\quad \frac{u(\partial)}{\alp},\quad x =\partial. 	
\end{cases}
\end{align*}
Since $u\in C(X_{\partial})$, 
we get $u|_X\in C_b(X)$ and hence from the above formula together with Corollary \ref{Feller-Resolvent} we conclude that $\check{R}_\alpha u$ is continuous at $x\in X$. For $x=\partial$, if it is isolated from $X$, then $\check{R}_\alpha u$ is continuous at $\partial$; if it is not isolated  from $X$, 
let $(x_n)_n\subset X$ be such that $x_n\to\partial$. Then since $u|_X\in C_b(X)$ we get by Corollary \ref{Feller-Resolvent}, once again, we have  $\check{R}_\alpha (u|_X)=\check{R}_{\alpha}(u|_X)
\in C_0(X)$ and $\check{R}_\alpha 1\in C_0(X)$ and hence
\[
\lim_{n\to\infty} \check{R}_\alpha (u|_X)(x_n) = \lim_{n\to\infty} \check{R}_\alpha 1(x_n) = 0,
\]
leading to $\lim_{n\to\infty} \check{R}_\alpha u(x_n) = \frac{u(\partial)}{\alp} =  \check{R}_\alpha u(\partial)$. We thus have proved $\check{R}_\alp(C(X_\partial))\subset C(X_\partial)$.\\
Set $F_\partial:= F\cup\{\partial\}$ and
\[
\sigma_{F_\partial}:= \inf\{t>0: \mathbb{X}_t\in F_\partial\},\quad  
P_{F_\partial} u(x) = \mathbb{E}_x\left[u(\mathbb{X}_{\sigma_{F_\partial}})\right], \quad x\in X_\partial,\quad u\in C_b(X_\partial).
\]
Let $u\in C(X_\partial)$. Then by the right-continuity of the sample paths $\check{\mathbb{X}}_t$ we obtain
\[
\lim_{\alp\to\infty} \alp\check{R}_\alp u(x) = \lim_{t\to 0}\check{\bf P}_tu(x) 
= \mathbb{E}_x\left[ u(\mathbb{X}_{\tau_0})\right],\quad x\in X_\partial.
\]
Let us show that $\mathbb{E}_x\left[ u(\mathbb{X}_{\tau_0})\right] = P_{F_\partial} u(x)$ for all $x\in X_\partial$.\\
Indeed, for $x\in X_\partial$ we have
\begin{align*}
\mathbb{E}_x\left[ u(\mathbb{X}_{\tau_0})\right] & = \mathbb{E}_x\left[ u(\mathbb{X}_{\sigma_F})\right] \\
& = \mathbb{E}_x\left[ u(\mathbb{X}_{\sigma_F}): \sigma_F>\sigma_{\{\partial\}}\right] + \mathbb{E}_x\left[ u(\mathbb{X}_{\sigma_F}): \sigma_F\leq\sigma_{\{\partial\}} \right]\\
& =  \mathbb{E}_x\left[ u(\mathbb{X}_\infty): \sigma_F>\sigma_{\{\partial\}}\right] + \mathbb{E}_x\left[ u(\mathbb{X}_{\sigma_F}): \sigma_F\leq\sigma_{\{\partial\}} \right]\\
& =  \mathbb{E}_x\left[ u(\mathbb{X}_{\sigma_{\{\partial\}}}): \sigma_F>\sigma_{\{\partial\}}\right] + \mathbb{E}_x\left[ u(\mathbb{X}_{\sigma_F}): \sigma_F\leq\sigma_{\{\partial\}} \right]\\
& =  \mathbb{E}_x\left[ u(\mathbb{X}_{\sigma_F\wedge\sigma_{\{\partial\}}}): \sigma_F>\sigma_{\{\partial\}}\right] + \mathbb{E}_x\left[u(\mathbb{X}_{\sigma_F\wedge\sigma_{\{\partial\}}}) : \sigma_F\leq\sigma_{\{\partial\}} \right]\\
& =  \mathbb{E}_x\left[ u(\mathbb{X}_{\sigma_F\wedge\sigma_{\{\partial\}}})\right]
= \mathbb{E}_x\left[u(\mathbb{X}_{\sigma_{F_\partial}})\right] = P_{F_\partial} u (x).
\end{align*} 
We have thus achieved 
\begin{equation}
\lim_{\alp\to\infty} \alp\check{R}_\alp u(x) = P_{F_\partial} u (x), \quad u\in C(X_\partial), \quad x\in X_\partial.
\label{Pointwise-Res-Hitting}
\end{equation}
As  $F_\partial$ is closed in $X_\partial$, using   \cite[Th\'eor\`eme IV, p.~425]{Lion} (where the set $F$ is denoted by  $A$) we get that $P_{F_\partial} u\in C(X_\partial)$ for all $u\in C(X_\partial)$ and the  point-wise  convergence from (\ref{Pointwise-Res-Hitting}) holds in fact  uniformly on $X_\partial$:
\begin{equation}
\lim_{\alp\to\infty}\sup_{x\in X_\partial} \left|\alp\check{R}_\alp u - P_{F_\partial} u \right| =0,\quad u\in C(X_\partial).
\label{uniform-extended}
\end{equation}
Besides we already proved that
\[
P_{F_\partial} u (x) = 
\begin{cases}
\quad P_Fu(x),\quad x\in X,\\
\quad u(\partial),\quad x=\partial,
\end{cases}
\quad u\in C(X_\partial).
\] 
Set $C_\infty(X_\partial):=\{u\in C(X_\partial): u(\partial) = 0\}$. Then, from the expression of $P_{F_\partial}$ we derive 
\begin{equation}
P_{F_\partial}(C_\infty(X_\partial)) \subset C_\infty(X_\partial).
\label{inclusion-1}
\end{equation}
Let $u\in C_0(X)$. Using our convention concerning extension of functions from $X$ to $X_\partial$ by setting $\tilde{u}(\partial) = 0$, we get $\tilde u\in C_\infty(X_\partial)$. Conversely, for any $u\in C_\infty(X_\partial)$  it holds $u|_X\in C_0(X)$, so that we can and  will identify $C_0(X)$ with $C_\infty(X_\partial)$. With this identification
in conjunction with the observation $P_F u = P_{F_\partial} \tilde{u}|_F$, the uniform convergence (\ref{uniform-extended}) and the inclusion (\ref{inclusion-1}) we get
\begin{equation}
P_F(C_0(X)) \subset C_0(X)\ \quad \text{and}\quad  \lim_{\alp\to \infty} \| \alp\check{R}_\alp u - P_F u \|_\infty = 0,\quad u\in C_0(X),
\label{inclusion-U-convergence}
\end{equation}
and the proof is finished.
\end{proof}

As a first consequence of Lemma~\ref{Feller-1} we get an improvement of Theorem \ref{right} in the following way:
\begin{coro}
The equivalent assertions from Theorem \ref{right} are equivalent to: 
\begin{enumerate}
\item[\rm(5)] $F$ is closed and $\check{R}_\alp$ is injective on $C_b(X)$ for each $\alp\geq 0$.
\end{enumerate}
\end{coro}
\begin{proof}
It suffices to prove  that $(4)$ from Theorem \ref{right} is equivalent to $(5)$.\\
$(5)\Rightarrow (4)$: Assume that $\check{R}_\alp$ is injective on $C_b(X)$ and let $u\in C_0(X)$ be such that $P_F u =0$. Then $\check{R}_\alp P_Fu = \check{R}_\alp u =0$ and hence $u=0$, by assumption.\\
Now let $u\in C_0(X)$; since $F$ is closed by assumption, we get from Lemma \ref{Feller-1} that the function $v:=P_Fu - u\in C_0(X)$; moreover $P_Fv =0$. Thus from the first step we get $v=0$ and hence $P_Fu = u$. This leads to $P_Fu = u$ for all $u\in C_0(X)$ and then for all indicator function $u=1_G$ of open set $G$ (by monotonous approximation) and then for all $u\in \mathscr{B}_b(X)$ (by monotone class argument), which is equivalent to injectivity  of $P_F$.\\
$(4)\Rightarrow (5)$: (4) is equivalent to (2) by Theorem~\ref{right}, in particular, we get the closedness of  $F$. The latter assertion is obvious with the help of Lemma \ref{sg-res-limit}.
\end{proof}
\begin{lem}
Assume that $F$ is closed. Then for all $\alp>0$ we have  $\overline\ran(\check{R}_\alp)=\overline\ran(G^\mu) = P_F(C_0(X))$.
\label{Range}
\end{lem}
\begin{proof}
The proof works as the proof of \cite[Lemma 2.5]{Duncan}.\\
From the fact that $\overline\ran(\check{R}_\alp)$ is independent of $\alp$ and that $\lim_{\alp\to 0} \check{R}_\alp= G^\mu$ in the operator norm, using the resolvent formula (\ref{Resolvent-Formula}) we deduce $\overline{\ran}(G^\mu)\subset \overline\ran(\check{R}_\alp)$. The reversed inclusion is obtained from (\ref{Resolvent-Formula}), as $\check{R}_\alpha = (1+\alpha G^\mu)^{-1} G^\mu= G^\mu (1+\alpha G^\mu)^{-1}$.\\
Let $u\in C_0(X)$. Then $P_F G^\mu u = G^\mu u$. Indeed, 
\begin{align*}
P_FG^\mu u(x) & = \mathbb{E}_x\left[  \mathbb{E}_{\mathbb{X}_{\sigma_F}}\left[ \int_0^\infty u(\mathbb{X}_t)\,d{\bf A}_t \right] \right]\\
& = \mathbb{E}_x\left[  \mathbb{E}_x\left[\left. \int_0^\infty u(\mathbb{X}_{t+\theta_{\sigma_F}})\,d{\bf A}_{t+\theta_{\sigma_F}}\,\right|\, {\mathscr{M}_{\sigma_F}}\right] \right]\\
& = \mathbb{E}_x\left[ \int_{\sigma_F}^\infty u(\mathbb{X}_t)\,d{\bf A}_t \right] = \mathbb{E}_x\left[ \int_0^\infty u(\mathbb{X}_t)\,d{\bf A}_t \right]\\
& = G^\mu u(x). 
\end{align*} 
The latter identity leads in conjunction with Lemma \ref{Feller-1} to $\overline\ran(G^\mu)\subset P_F(C_0(X))$.\\
Conversely, let $u\in C_0(X)$. By Lemma \ref{Feller-1} we get $\lim_{\alpha\to\infty}\alp\check{R}_\alpha u =P_F u$ uniformly. As $\ran(\check{R}_\alpha)$ is independent of $\alpha$ we thereby get $P_F u\in \overline\ran(\check{R}_\alpha)=\overline\ran(G^\mu)$. 
\end{proof}
\begin{theo}[A substitute to Feller property, the transient case] Assume that $F$ is closed. Then
\begin{enumerate}
\item[\rm(1)] The space $P_F(C_0(X))$ is  $\check{R}_\alpha$-invariant for all $\alp>0$ and $\lim_{\alp\to\infty} \|\alp\check{R}_\alp u - u\|_\infty = 0$ for all $u\in P_F(C_0(X))$. 
\item[\rm(2)] The space $P_F(C_0(X))$ is $\check{\bf P}_t$-invariant for all $t>0$ and $\lim_{t\to 0}\| \check{\bf P}_t u - u \|_\infty = 0$ for all $u\in P_F(C_0(X))$.
\item[\rm(3)] For all $t\geq 0$ the semigroup $\check{\bf P}_t$ maps continuously $C_0(X)$ into itself. Moreover, it is continuous, i.e., for all $u\in C_0(X)$, the map: $t\in[0,\infty)\to C_0(X),\; t\mapsto \check{\bf P}_t u$ is continuous.
\end{enumerate}
\label{Substitute-Feller}
\end{theo}
\begin{proof}
Let $\alp,\beta,t>0$.  By Lemma  \ref{Range} we have  $\overline{\ran}(\check{R}_\beta)=P_F(C_0(X))$ and hence using the property  $\check{R}_\alpha\check{R}_\beta = \check{R}_\beta \check{R}_\alpha$ we get that it is $\check{R}_\alpha$-invariant. Markov property implies that  $\check{R}_\alpha$ commutes with $\check{\bf P}_t$ and hence $P_F(C_0(X))$ is also invariant for $\check{\bf P}_t$.\\ 
Let us consider the Banach space $P_F(C_0(X))$ equipped with the topology of uniform convergence. Owing to Proposition \ref{First-Prop-Res}-(1) in conjunction with Lemma \ref{Feller-1}-(2) we conclude that the family  $\check{R}_\alpha:P_F(C_0(X))\to P_F(C_0(X))$ is a strongly continuous sub-Markovian resolvents family. An appeal to Hille--Yosida theorem yields that $\check{\bf P}_t$ preserves  $P_F(C_0(X))$ and $\check{\bf P}_t:P_F(C_0(X))\to P_F(C_0(X))$ is a strongly continuos sub-Markovain semigroup. We have thus proved (1)-(2).\\
For the remaining part of the proof, let $u\in C_0(X), t_0\geq 0$. As $\check{\mathbb{X}}_t\in F$, $\mathbb{P}_x$-a.s. and $P_Fu = u$ on $F$, we get  $\check{\bf P}_t P_F u = \check{\bf P}_t u\in C_0(X)$ for all $t\geq 0$, by the second assertion and Lemma \ref{Feller-1}.\\
By (2) we get  $\lim_{t\to t_0}\| \check{\bf P}_t P_F u - \check{\bf P}_{t_0} P_F u \|_\infty =0$. 
Using the identity $\check{\bf P}_t P_F u = \check{\bf P}_t u$, once again, we obtain $0=\lim_{t\to t_0}\| \check{\bf P}_t P_F u - \check{\bf P}_{t_0} P_F u \|_\infty = \lim_{t\to t_0}\| \check{\bf P}_t  u - \check{\bf P}_{t_0}  u \|_\infty$.\\
\end{proof}
Theorem \ref{Substitute-Feller} can be interpreted as a substitute to Feller property. Indeed, the theorem states that $\check{\bf P}_t$ is continuous on  $C_0(X)$, but is not a $\mathbf{C}_0$-semigroup on $C_0(X)$. The following special case, which is an immediate consequence of the theorem, supports our interpretation.
\begin{coro}
Assume that $F=X$. Then we have the following:
\begin{enumerate}
\item[\rm(1)] The resolvent $(\check{R}_\alpha)_{\alpha>0}$ has the Feller property on $C_0(X)$.	
\item[\rm(2)] The semi-group $(\check{{\bf P}}_t)_{t\geq0}$ has the Feller property on $C_0(X)$.
\end{enumerate}
\label{Feller-Full}
\end{coro}
By the end of this section we mention that though the family $\check{R}_\alp$ is a (pseudo)-resolvent, i.e., a family $(\check{R}_\alp)_{\alpha>0}$
of continuous linear contractions which satisfies the resolvent identity but which might be non-injective and the semigroup $\check{{\bf P}}_t$ is degenerate, i.e., a family $(T_t)_{t\geq0}$  of continuous linear contractions which satisfies $T_tT_s = T_{t+s}$, the map $t\in [0,\infty)\mapsto T_t u$ is continuous  for all $u$ but  $T_0$ is not the identity map,  nevertheless they still generate a closed linear relation in the following sense (see \cite[Theorem 2.1]{Baskakov-1}, \cite[Theorem 5.2]{Arendt-23}): Let $\check{L}_\rel$ be the closed linear relation defined by the linear subspace of $C_0(X)\times C_0(X)$
\[
\check{L}_\rel:= \{ (u, \check{R}_\alp^{-1}\{u\} - \alp u): u\in C_0(X) \}. 
\]
Then every $\alp>0$ is in the resolvent set of $-\check{L}_\rel$ moreover 
\[
(\check{L}_\rel +\alp)^{-1} = \check{R}_\alp\quad \text{and}\quad e^{-t \check{L}_\rel} = \check{{\bf P}}_t,\; t\geq0.
\]
We will say that $-\check{L}_\rel$ is the generator of $\check{R}_\alp$ or  $\check{{\bf P}}_t$.\\
For any $u\in\dom(\check{L}_\rel) = P_F(C_0(X))$ we denote by $\check{L}_\rel u$ the set $\{ (u, \check{R}_\alp^{-1}\{u\} - \alp u) \}$. 
\begin{prop}
Let $-\check{L}_\rel$ be the generator of $(\check{R}_\alp)_\alp$. Then we have the following:
\begin{enumerate}
\item[\rm(1)] The linear relation  $-\check{L}_\rel$ is $m$-dissipative.
\item[\rm(2)] The linear relation  $-\check{L}_\rel$ generates  an integrated semigroup, $\check{\bf S}_t: C_0(X)\to C_0(X)$, $t\geq 0$ via the equation
\begin{equation}
\check{R}_\alp = \alp\int_0^\infty e^{-\alp t} \check{\bf S}_t\,dt.
\label{1-integrated}
\end{equation} 
\item[\rm(3)] It holds
\begin{equation}
	\check{\bf S}_t = \int_0^t \check{\bf P}_s\,ds,\quad t\geq 0.
	\label{integrated-smgp}
\end{equation}
\end{enumerate}
Here both integrals in the right hand side of 
\eqref{1-integrated} and \eqref{integrated-smgp}, respectively are
 the strong Riemann integrals.
\end{prop}
\begin{proof}
The first two assertions follow  from Proposition \ref{First-Prop-Res} and \cite[Proposition 4.3, Theorem 5.9]{Arendt-23}.\\
The last equation follows from continuity of $t\mapsto\check{\bf P}_t u$ for all $u\in C_0(X)$ in conjunction with formula (\ref{1-integrated}) and integration by parts, which completes the proof.
\end{proof}
By the end of this section, let us explain how to extend the results of this section to general (not necessary transient) standard processes.\\
Let $\mathbb{X}$ be  a standard process,  $q>0$ and $\mathbb{X}^q$ the subprocess associated to the semigroup ${\bf Q}_t^q:=e^{-qt}{\bf P}_t$. Then $\mathbb{X}^q$ is transient. Let $\mu$ be a positive Borel smooth measure on subsets of $X$ with fine support $F$ and ${\bf A}$ its PCAF (see \cite[Th\'eor\`eme VI. 1, p.~524]{Revuz}). Set $\check{\mathbb{X}}^q$ the time change of $\mathbb{X}^q$ by means of ${\bf A}$. Then the semigroup and the resolvent of $\check{\mathbb{X}}^q$ are respectively:
\[
\check{R}_\alpha^q u(x):= \mathbb{E}_{x} \left[ \int_0^\infty  e^{-\alpha t-q\tau_t}u(\mathbb{X}_{\tau_t})\,dt \right],\ \check{\bf Q}_t^q u(x):=\mathbb{E}_{x} \left[ e^{-q\tau_t}u(\mathbb{X}_{\tau_t}) \right],\; x\in X,\; u\in\mathscr{B}_b(X),
\]
and for any $u\in C_b(X)$ it holds
\[
\lim_{t\downarrow 0} \check{\bf Q}_t^q u(x) = \lim_{\alpha\to\infty} \alpha \check{R}_\alpha^q u(x) = \mathbb{E}_x( e^{-q\sigma_F} u(\mathbb{X}_{\sigma_F} ))=: P_F^q u(x),\; x\in X.
\]
Assume further that for some (and hence every) $\alpha>0$ the measure $\mu$ is  $G_\alpha$-Kato and $G_\alpha(\cdot,y)$ is l.s.c for any $y\in X$. Using Corollary \ref{Feller-Resolvent} we obtain 
\begin{equation}
 \check{R}_\alpha^q:\mathscr{B}_b(X)\to C_0(X)\ \text {continuously for any}\  \alpha,q>0.
\label{smooth-a-resolvent} 
\end{equation}
Assume that $F$ is closed. By Lemma \ref{Feller-1} we get for all $q>0$, $u\in C_0(X)$
\begin{equation}
 P_F^q(C_0(X))\subset C_0(X)\ \text{and}\ \lim_{\alpha\to\infty} \alpha \check{R}_\alpha^q u =  P_F^q u\ \text{uniformly}.
\label{Feller-2}
\end{equation} 
 Moreover Lemma \ref{Range} yields 
\begin{equation}
 \overline{\ran}( \check{R}_\alpha^q ) = \overline{\ran}(G_q^\mu) = P_F^q (C_0(X))\ \text{ for all}\ q>0.
 \label{Range-2}
\end{equation} 
Suppose  for a while that $X$ is  compact and let $u\in C(X)$. Then for any $\alpha>0$, we have
\begin{align}
\big|\check{R}_\alpha u(x) -  \check{R}_\alpha^q u(x) \big|& =   
\left|\mathbb{E}_x \left[ \int_0^\infty (1 - e^{-q\tau_t})e^{-\alpha t} u(\mathbb{X}_{\tau_t})\,dt\right] \right|\\
&\leq \|u\|_\infty \left( \frac{1}{\alpha} - \mathbb{E}_x \left[ \int_0^\infty e^{-\alpha t - q\tau_t} \,dt \right] \right).
\label{resolv-control}
\end{align}
By dominated convergence theorem we get that the decreasing, w.r.t $q$, family of functions
\[ 
g_q(x) := \frac{1}{\alpha} - \mathbb{E}_x \left[ \int_0^\infty e^{-\alpha t - q\tau_t} \,dt \right]
= \frac{1}{\alpha} - \check{R}_\alpha^q 1(x)
\]
decrease point-wisely to $0$ as $q$ decreases to $0$ for any $\alpha>0$. By (\ref{smooth-a-resolvent}), $(g_q)_{q>0}\subset  C(X)$, so that by Dini's lemma $\lim_{q\to 0} g_q= 0$ uniformly and hence by (\ref{resolv-control})
\begin{equation}
 \lim_{q\to 0} \check{R}_\alpha^q u = \check{R}_\alpha u\quad \text{uniformly}.
 \label{q-Res-Convergence}
\end{equation} 
Since we already know that $\check{R}_\alpha^q u\in C(X)$ for all $q>0$, we conclude that $\check{R}_\alpha u\in C(X)$ and hence $\check{R}_\alpha u : C(X)\to C(X)$, continuously. Now an appeal to Lion's theorem \cite[Th\'eor\`eme IV, p.425]{Lion} leads to
\begin{equation}
\lim_{\alpha\to\infty} \alpha \check{R}_\alpha u = P_F u\quad \text {uniformly}.
\label{hitting-limit}
\end{equation}
Using the first part of (\ref{Feller-2}), we obtain
\begin{equation}
\lim_{q\to 0} P_F^q u = P_F u\quad \text{ uniformly},
\label{-q-hitiing-Convergence}
\end{equation}
so that  $P_F(C(X))\subset C(X)$ and $\lim_{q\to 0} \check{{\bf Q}}_t^q u = \check{{\bf P}}_t u $ uniformly.\\
Let us show that $\overline{\ran} (\check{R}_\alpha) = P_F(C(X))$. From (\ref{hitting-limit}), we derive $P_F(C(X)) \subset  \overline{\ran}(\check{R}_\alpha)$. By Proposition \ref{ResolventFormula}, we have
\[
\check{R}_\alpha^q  = (1+\alpha G_q^\mu)^{-1} G_q^\mu = G_q^\mu(1+\alpha G_q^\mu)^{-1}, 
\] 
which leads in conjunction with the formula  $P_F^q G_q^\mu  = G_q^\mu$ (see the proof of Lemma \ref{Range}) to $ P_F^q \check{R}_\alpha^q u = \check{R}_\alpha^q u, q>0$, $u\in C(X)$. Thus passing to the limit  $q\to 0$ and using (\ref{q-Res-Convergence})-(\ref{-q-hitiing-Convergence}) we obtain $P_F\check{R}_\alpha u = \check{R}_\alpha u$. Recalling that $\check{R}_\alpha(C(X))\subset C(X)$ we thereby obtain $\overline{\ran}(\check{R}_\alpha) \subset P_F(C(X))$. Hence  $\overline{\ran} (\check{R}_\alpha) = P_F(C(X))$.\\
We have thus proved that if $\mathbb{X}$ is a standard process (not necessarily transient), $X$ is compact, $\mu$  is $G_\alpha$-Kato, $G_\alpha(\cdot,y), y\in X$ is l.s.c for some $\alpha>0$ and $F$ is closed, then the assertions of  Lemma \ref{Feller-1} hold true and $\overline{\ran} (\check{R}_\alpha) = P_F(C(X))$. Now mimicking the proof of Theorem \ref{Substitute-Feller} we conclude that the conclusions of the latter  still  hold true in this situation.\\
Now if $X$ is a general set (not necessarily compact), going over to the compactification $X_\partial$ of $X$, extending the resolvent and the semigroup to $X_\partial$ as it was done in the proof of Lemma \ref{Feller-1} and using the already proved results for the compact case we obtain:
\begin{theo}[A substitute to Feller property, the general case] Let $\mathbb{X}$ be a standard process. Assume that $\mu$ is $G_\alpha$-Kato, $G_\alpha(\cdot,y), y\in X$ is l.s.c for some $\alpha>0$ and that the fine support $F$ of $\mu$ is closed. Then, the assertions of Theorem \ref{Substitute-Feller} still hold true.
\label{Substitute-Feller-general}
\end{theo}

\section{Convergence results}
In this section we fix a sequence of measures  $(\mu_n)_{n\in\N_\infty}\subset\mathscr{K}(X)$ such that 
\begin{align}
	\mu_n\to\mu_\infty\quad \text{vaguely}\quad\text{and}\quad	G^{\mu_n} 1\to G^{\mu_\infty} 1\quad \text{uniformly}.
	\label{PotentialConvergence}		
\end{align}
Let $F_n$, $F$ be the fine supports of $\mu_n$, $\mu_\infty$, respectively.
\begin{rk}
{\rm
In \cite[Theorem 3.2]{BenAmor-24} it is proved that if the measures $\mu_n$ are finite and converge weakly monotonically, i.e. the $\mu_n$'s either increase or decrease weakly to $\mu_\infty$, then the condition $G^{\mu_n} 1\to G^{\mu_\infty} 1$ uniformly is automatically fulfilled.
}
\end{rk}
\begin{theo}
Let $(\mu_n)_{n\in\N_\infty}\subset\mathscr{K}(X)$ satisfy assumption (\ref{PotentialConvergence}). Then for all $u\in {C}_0(X)$, it holds
\begin{enumerate} 
\item[\rm(1)] 	$\lim_{n\to\infty}\|G^{\mu_n} u - G^{\mu_\infty} u\|_\infty =0$.
\item[\rm(2)] $\lim_{n\to\infty}\|\check{R}_\alp^n u - \check{R}_\alp^\infty  u\|_\infty =0$ for all $\alp\geq 0$.
\end{enumerate}	
\label{UniformConvergencePotOPerator}
\end{theo}
\begin{proof}
The second assertion follows from the first one in conjunction with formula (\ref{Resolvent-Formula}). Let us prove the first assertion.\\	
Since ${C}_c(X)$ is dense in ${C}_0(X)$, it suffices to prove uniform convergence on ${C}_c(X)$.\\
Let $u\in {C}_c(X)$, without loss of generality we may and will assume that $0\leq u\leq 1$.\\
Since $G(\cdot,\cdot)$ is lower semi-continuous and $G\not\equiv\infty$, there is a sequence $(H_k(\cdot,\cdot))_k\subset {C}_b(X\times X)$ such that  $H_k(x,y)\uparrow G(x,y)$ for all $x,y\in X$ (see \cite[Theorem 2.1.3, p.~26]{Ransford}). Thus, for each fixed $n\in\N_\infty$ we get $\int_X  H_k(\cdot,y)u(y)\,d\mu_n(y)\uparrow \int_X  G(\cdot,y)u(y)\,d\mu_n$. By Proposition \ref{PotentialOperator} we conclude that the functions  $\int_X  H_k(\cdot,y)u(y)\,d\mu_n(y)$ are in $C_0(X)$  and accordingly
\[
	\lim_{k\to\infty}\int_X  H_k(\cdot,y)u(y)\,d\mu_n(y) =\int_X  G(\cdot,y)u(y)\,d\mu_n(y)\quad \text{uniformly},
\]
in view of Dini's lemma on $C(X_\partial)$.\\
As $H_k(x,\cdot)u\in {C}_c(X)$ we get, by vague convergence of the $\mu_n$'s towards $\mu_\infty$, for each fixed integer $k$ it holds
	\[
	\int_X  H_k(\cdot,y)u(y)\,d\mu_n(y)\to \int_X  H_k(\cdot,y)u(y)\,d\mu_\infty(y)\quad  \text{point-wisely}\quad\text{as}\quad n\to\infty.
	\]
Moreover, from continuity of $H_k$ we get that the family $\int_X  H_k(\cdot,y)u(y)\,d\mu_n(y),n\in\N$ is equicontinuous on each compact set and is uniformly bounded, for fixed $k$. Thus
\begin{equation}
		\lim_{n\to\infty} \int_X  H_k(\cdot,y)u(y)\,d\mu_n(y) = \int_X  H_k(\cdot,y)u(y)\,d\mu_\infty(y)\quad \text{locally uniformly}.
		\label{LocUniform}
\end{equation}
Now  let $\epsilon>0$ and  $K\subset X$ a compact set such that $1_{K^c}G^{\mu_\infty}1<\epsilon$. Set
\[
	v_{k,n}(x):= \int_X  H_k(x,y)u(y)\,d\mu_n(y), \quad x\in X.
\]
We write $v_{k,n} = 1_K v_{k,n} + 1_{K^c}v_{k,n}$. On the one hand we have, by the convergence assumption $ 1_{K^c} v_{k,n}\leq 1_{K^c}G^{\mu_n} 1\to 1_{K^c}G^{\mu_\infty} 1<\epsilon$, $1_{K^c} v_{k,\infty} \leq 1_{K^c} G^{\mu_\infty}1<\epsilon$. Thus
\[
	\sup_k \left( 1_{K^c} v_{k,n} \right)< 2\epsilon\quad \text{for large}\ n, \quad 1_{K^c} v_{k,\infty}<\epsilon.
\]  
On the other hand we have by (\ref{LocUniform}), for any fixed $k$,
\[
	\lim_{n\to\infty} 1_K v_{k,n} =  1_K v_{k,\infty}\quad \text{uniformly}.
\]
Finally let $\epsilon>0$ and $k$ be such that $\|G^{\mu_\infty} u  - v_{k,\infty}\|_\infty <\epsilon$. Then
\begin{align*}
		\|G^{\mu_\infty} u - G^{\mu_n} u \|_\infty& \leq \|G^{\mu_\infty} u  - v_{k,\infty}\|_\infty + \|v_{k,\infty} - v_{k,n}\|_\infty + \|v_{k,n} - G^{\mu_n}1\|_\infty\\
		&\leq \epsilon + \|v_{k,\infty} - v_{k,n}\|_\infty + \|v_{k,n} - G^{\mu_n}1\|_\infty.
\end{align*}
We thereby derive
\begin{align*}
		\lim_{n\to\infty }\|G^{\mu_\infty} u - G^{\mu_n} u \|_\infty \leq \epsilon + \lim_{n\to\infty}\|v_{k,\infty} - v_{k,n}\|_\infty 
		+ \lim_{n\to\infty}\|v_{k,n} - G^{\mu_n}1\|_\infty = 4\epsilon.
\end{align*}
As $\epsilon$ is arbitrary we obtain $\lim_{n\to\infty }\|G^{\mu_\infty} u - G^{\mu_n} u \|_\infty =0$, which ends the proof.
\end{proof}
Now using the latter theorem together with \cite[Theorem 5.17]{Arendt-23} we obtain
\begin{coro}
Let $\check{\bf S}_t^n, n\in\N_\infty$ be the integrated semigroup associated with $\check{R}_\alp^n$ on the space $C_0(X)$. Then for all $u\in C_0(X)$ it holds  $\lim_{n\to\infty} \|\check{\bf S}_t^n u - \check{\bf S}_t^\infty u\|_\infty =0$ locally uniformly in $t\in[0,\infty)$. 
\label{Conv-Integrated} 
\end{coro}
We come now to convergence of semigroups.
\begin{theo}\label{thm:ConvergenceSemiGroup}
\begin{enumerate}
\item[\rm(1)] Let $u\in\overline\ran(\check{R}_\alp^\infty)$. Then $\lim_{n\to\infty} \|\check{\bf P}_t^n u - \check{\bf P}_t^\infty u\|_\infty =0$  for all  $t\geq 0$. Moreover, the convergence is uniform in $t$ on compact intervals of $[0,\infty)$. 
\item[\rm(2)] Assume that $F$ is closed. Let $u\in C_0(X)$. Then $\lim_{n\to\infty} \|\check{\bf P}_t^n P_Fu - \check{\bf P}_t^\infty u\|_\infty=0$ for all $t\geq 0$. Moreover, the convergence is uniform in $t$ on compact intervals of $[0,\infty)$. 
\end{enumerate}
\label{Convergence-Sg}
\end{theo}
Though the proof can be obtained by mimicking the proof of \cite[Theorem 1.1]{Duncan} (except uniform local convergence in $t$), we give here a shorter one.
\begin{proof}
(1): Let $n\in\N_\infty$. By Hille--Yosida theorem the semigroup $\check{\bf P}_t^n$ is strongly continuous on $\overline\ran(\check{R}_\alp^n)$ and hence  $\check{\bf P}_t^n u$ is strongly differentiable in $t$ for all $u$ in the domain of the generator of  $\check{\bf P}_t^n$ restricted to $\overline\ran(\check{R}_\alp^n)$. Thus for all $u\in C_0(X)$ the function  $\check{\bf P}_t^n \check{R}_\alp^n u$ is strongly differentiable in $t$. Let $u\in C_0(X)$. Following Kato's computation (see \cite[p.~503, identity (2.27)]{Kato}) we get
\[
\check{R}_\alp^n (\check{\bf P}_t^\infty - \check{\bf P}_t^n ) \check{R}_\alp^\infty u = 
\int_0^t  \check{\bf P}_{t-s}^n (\check{R}_\alp^\infty -   \check{R}_\alp^n) \check{\bf P}_s^\infty u\,ds.
\]
Owing to strong convergence of the resolvents in conjunction with dominated convergence theorem we achieve $\lim_{n\to\infty}\|\check{R}_\alp^n (\check{\bf P}_t^\infty - \check{\bf P}_t^n ) \check{R}_\alp^\infty u\|_\infty = 0$ and the convergence holds locally uniformly in $t$. Thus for all $u\in \overline\ran(\check{R}_\alp^\infty)$ we get $\lim_{n\to\infty}\|\check{R}_\alp^n (\check{\bf P}_t^\infty - \check{\bf P}_t^n ) u\|_\infty = 0$.\\
Recalling that $\lim_{n\to\infty}\|\check{R}_\alp^n \check{\bf P}_t^\infty u -  \check{R}_\alp^\infty \check{\bf P}_t^\infty u\|_\infty = 0$ we obtain 
\[
\lim_{n\to\infty}\|\check{R}_\alp^n \check{\bf P}_t^n u - \check{\bf P}_t^\infty \check{R}_\alp^\infty u\|_\infty = 0.
\]
Observing that $\check{R}_\alp^n \check{\bf P}_t^n u  = \check{\bf P}_t^n \check{R}_\alp^n u$ and 
$\lim_{n\to\infty}\|\check{\bf P}_t^n\check{R}_\alp^n  u - \check{\bf P}_t^n \check{R}_\alp^\infty u\|_\infty = 0$ we conclude that 
$\lim_{n\to\infty}\|\check{\bf P}_t^n \check{R}_\alp^\infty u  -  \check{\bf P}_t^\infty \check{R}_\alp^\infty u \|_\infty = 0$ locally uniformly in $t\in[0,\infty)$ and thereby $\lim_{n\to\infty}\|\check{\bf P}_t^n  u  -  \check{\bf P}_t^\infty u \|_\infty = 0$ for all $u\in \overline\ran(\check{R}_\alp^\infty)$  locally uniformly in $t\in[0,\infty)$ and the first assertion is proved.\\
(2): If furthermore $F$ is closed, by Lemma \ref{Range} we have $\overline\ran(\check{R}_\alp^\infty) = P_F(C_0(X))$. Observing that $ \check{{\bf P}}_t^\infty P_F u = \check{{\bf P}}_t^\infty u$ and applying the first assertion of the theorem we get the result.
\end{proof}
In some special cases we have very interesting consequences of the latter theorem.
\begin{prop}\label{prop:ConvergenceFeller}
Assume that $F_n,\,F$ are closed and   $F_n\subset F$ for every $n$. Then for all $u\in C_0(X)$ we have $\lim_{n\to\infty} \|\check{\bf P}_t^n u - \check{\bf P}_t^\infty u\|_\infty = 0$ locally uniformly in $t\in[0,\infty)$.
\label{SpecialCase}
\end{prop}
\begin{proof}
	 Follows from 	$\check{{\bf P}}_t^n P_F u = \check{{\bf P}}_t^n  u$ in conjunction with Theorem \ref{Convergence-Sg}.
\end{proof}
\begin{coro}
Assume that $F_n,\,F$ are closed and $\mu_n\uparrow\mu_\infty$ vaguely. Then for all $u\in C_0(X)$ it holds  $\lim_{n\to\infty} \|\check{\bf P}_t^n u - \check{\bf P}_t^\infty u \|_\infty = 0$ locally  uniformly in $t\in[0,\infty)$. 	
\label{Monotone-Conv}
\end{coro}
\begin{proof}
The result follows from Proposition \ref{SpecialCase}-(1).
\end{proof}
\begin{coro}
Assume that $F=X$. Then for all $u\in C_0(X)$ it holds  $\lim_{n\to\infty} \|\check{\bf P}_t^n u - \check{\bf P}_t^\infty u \|_\infty =0$ locally  uniformly in $t\in[0,\infty)$. 
\label{Full-Support}
\end{coro}
\begin{proof}
Since $F=X$ we get $P_F u=u$, which in conjunction with Theorem \ref{Convergence-Sg}{\color{red}{-(2)}} yields the  assertion.
\end{proof}
\section{Applications}
In this section we still fix a sequence $(\mu_n)_{n\in\N_\infty}\subset\mathscr{K}(X)$ such that  $\mu_n\to\mu_\infty$ vaguely and $G^{\mu_n}1 \to G^{\mu_\infty}1$ uniformly.\\
The first application concerns convergence of solutions of evolution equations.  Let $(v_n)_{n\in\N_\infty}$ $\subset C_0(X)$. According to \cite[Corollary 5.16]{Arendt-23},  $\check{\bf S}_t^n v_n$ is the unique classical solution of the evolution equation
\begin{equation}
\begin{cases} 
	- \frac{\partial u_n(t)}{\partial t} \in \check{L}_\rel^n u_n(t) -  v_n,\\
	u_n(0) =0,
\end{cases}
\label{Evolution}
\end{equation}
where the relation  $\check{L}_\rel^n  -  v_n$ is defined by (see \cite[p.6]{Arendt-23}):
\[
\check{L}_\rel^n  -  v_n =\{ (x_1, x_2 - v_n), (x_1,x_2)\in \check{L}_\rel^n \}
\]
and $u_n(t)$ being a classical solution of (\ref{Evolution}) means that
\[ u_n\in C^1((0,\infty), C_0(X))\cap C([0,\infty), C_0(X))\ \text{and}\  u_n\ \text{satisfies}\  (\ref{Evolution}).
\]
The latter is in turn equivalent to $u_n\in C^1((0,\infty), C_0(X))\cap C([0,\infty), C_0(X)), u_n(0)=0$ and
\[
\left( -\int_0^t u_n(s)\,ds, u_n(t) - tv_n \right)\in\check{L}_\rel^n\quad \text{for all}\quad t\geq 0. 
\]
\begin{theo}
Assume that $\lim_{n\to\infty} \|v_n - v_\infty\|_\infty =0$, then $\lim_{n\to\infty} \|\check{\bf S}_t^n v_n - \check{\bf S}_t^\infty v_\infty \|_\infty =0$ locally uniformly in $t\in [0,\infty)$
\label{Conv-Evolution-1}
\end{theo} 
\begin{proof}
From formula (\ref{integrated-smgp}) and since $\check{{\bf P}}_t^n$ are sub-Markovian we get $\|\check{\bf S}_t^n\| \leq t\| \check{{\bf P}}_t^n \|\leq t$. Hence using Corollary \ref{Conv-Integrated} and the triangle inequality we obtain the result.
\end{proof}
%
Let us consider now the heat equation. Assume further that $\check{\bf P}_t^n$ are holomorphic for all $n\in\N_\infty$.  Let $(v_n)_{n\in\N_\infty}\subset C_0(X)$. Then $\check{{\bf P}}_t v_n$ solves the heat equation 
\begin{equation}
	\begin{cases} 
		- \frac{\partial u_n(t)}{\partial t} \in \check{L}_\rel^n u_n(t),\\
		u_n(0) =v_n.
	\end{cases}
	\label{Heat-Eq-1}
\end{equation}
\begin{theo}
Assume that $\check{\bf P}_t^n$ are holomorphic for all $n\in\N_\infty$. Let $u_n(t), v_n$ be as before and $F_n,F$ be the respective fine supports of $\mu_n,\mu_\infty$. Assume that $\lim_{n\to\infty} \|v_n - v_\infty\|_\infty =0$, $F_n,F$ are closed and that either $F_n\subset F$ for all $n$ or $\mu_n\uparrow\mu_\infty$ vaguely. Then $\lim_{n\to\infty} \|u_n(t) - u_\infty(t)\|_\infty =0$ locally uniformly in $t\in [0,\infty)$.
\label{Conv-Heat}
\end{theo}
\begin{proof}
The proof follows from either Proposition \ref{SpecialCase} or Corollary \ref{Monotone-Conv}, with the help of the triangle inequality.
\end{proof}
\begin{rk}
{\rm 
We emphasize that  Theorems \ref{Conv-Evolution-1} and \ref{Conv-Heat} are new and relevant even in the case $F_n = F= X$ for all $n\in\N$. In this situation the linear relations $\check{L}_\rel^n$ become linear operators, the evolution equations (\ref{Evolution}) become  ordinary differential equations and $S_n(t)v_n$ are their respective classical solutions. Whereas equations (\ref{Heat-Eq-1}) become the usual heat equations,  $\check{{\bf P}}_t v_n$ are their respective classical solutions and the semigroups $\check{{\bf P}}_t^n, t>0$ become Feller semigroups, by Corollary \ref{Feller-Full}. Accordingly, Theorems \ref{Conv-Evolution-1} and \ref{Conv-Heat} express the fact that specific changes of the reference measure lead to continuous dependence of solution of the evolution equations (\ref{Evolution}) and of the heat equations (\ref{Heat-Eq-1}). 
}
\end{rk}
Let us now give some applications concerning the process $\check{\mathbb{X}}$ with state space $F$. We start by establishing some properties of the related semigroup. Let us emphasize that in the subsequent lines till Theorem \ref{Feller-restricted} we do not assume $\mathbb{X}$ to be transient.\\
Owing to the fact that for any $t\geq 0$ we have $\check{\mathbb{X}}_t\in F$, $\mathbb{P}_x$-a.s. for any $x\in X$, we conclude  that for any $t\geq 0$ the function $\check{\bf P}_tu $ depends solely on the values of $u$ on $F$:
\[
\check{\bf P}_tu(x) = \check{\bf P}_t(u|_{F})(x),\quad x\in X, \quad u\in\mathscr{B}_b(X), \quad t\geq 0. 
\]	
This leads us to reconsider the semigroup $\check{\bf P}_t$ on the space $\mathscr{B}_b(F)$.\\ From now on we assume that $F$ is closed. For any $t\geq 0$ we let
\[
\check{T}_t : \mathscr{B}_b(F)\to \mathscr{B}_b(F),\quad u\mapsto \mathbb{E}_{\cdot}\left[u(\mathbb{X}_{\tau_t})\right].
\]
It is easy to verify that $(\check{T}_t)_{t\geq 0}$ is a contractions semigroup. Moreover for any $u\in\mathscr{B}_b(F)$ and any Borel bounded extension $\tilde u\in\mathscr{B}_b(X)$ of $u$ we have 
\[
\check{T}_t u = (\check{\bf P}_t \tilde u)|_F.
\]
Let $(V_\alp)_{\alp\geq 0}$ be the resolvent of $\check{T}_t u$:
\[
V_\alp: \mathscr{B}_b(F)\to \mathscr{B}_b(F),\quad  V_\alp u(x) = \int_0^\infty e^{-\alp t}\check{T}_t u(x)\,dt,\quad u\in\mathscr{B}_b(F),\quad x\in F.
\]
It is easy to check that for all $u\in\mathscr{B}_b(F)$ one has $V_\alp u = \check{R}_\alp \tilde{u}|_F$ where $\tilde u$ is any Borel extension of $u$.
\begin{theo}\label{thm:Feller}
Let $\mathbb{X}$ be a (not necessarily transient) standard process. Suppose that for some $\alpha>0$, $\mu$ is $G_\alpha$-Kato, $G_\alpha(\cdot,y), y\in X$ is l.s.c and let $F$ be the fine support of $\mu$. Assume that $F$ is closed. Then
\begin{enumerate}
 \item[\rm(1)] The semigroup  $\check{T}_t$ maps $C_0(F)$ into itself and $\check{T}_t: C_0(F)\to C_0(F)$ has Feller property.
 \item[\rm(2)] The resolvent  $V_\alp$ maps $C_0(F)$ into itself and $V_\alp: C_0(F)\to C_0(F)$ has Feller property.
 \end{enumerate}
\label{Feller-restricted}
\end{theo}		
\begin{proof}
It suffices to prove the first assertion.\\
Let $u\in C_0(F)$. We claim that there is   $\tilde u\in C_0(X)$ such that  $\tilde{u}|_F = u$. Indeed, let $u\in C_0(F)$. Then the function $v:X_\partial\to\R$ such that $v=u$ on $F$ and $v(\partial) = 0$ satisfies $v\in C(F_\partial)$. By Tietze extension theorem there is $\tilde v\in C(X_\partial)$ such that   $\tilde{v}|_{F_\partial} = v$. As $\tilde{v}(\partial) = v(\partial) = 0$, we get $\tilde u:= \tilde{v}|_X\in C_0(X)$ and $\tilde{u}|_F = u$ and the claim is proved.\\ 
Thus using  Theorem \ref{Substitute-Feller-general} we get $\check{{\bf P}}_t \tilde u\in C_0(X)$. Closedness of $F$ yields $\check{T}_t u = (\check{\bf P}_t \tilde u)|_F\in C_0(F)$. Hence 	$\check{T}_t$ maps $C_0(F)$ into itself.\\
Besides using Theorem \ref{Substitute-Feller-general} once again we achieve
\begin{align*}
\sup_{x\in F} |  \check{T}_t u(x) - u(x) | = \sup_{x\in F} | \check{\bf P}_t \tilde u|_F(x) - \tilde u|_F(x) |
\leq \| \check{\bf P}_t \tilde u - P_F\tilde u\|_\infty \to 0\quad \text{as}\quad t\to 0.
\end{align*}
Finally it is obvious that $\check{T}_t: C_0(F)\to C_0(F)$ is a sub-Markovian semigroup.
\end{proof}	
\begin{coro}\label{cor:Feller}
Let $\Omega$ be a bounded Lipschitz domain of $\R^d$ with boundary $F$ and $\mu$ the $(d-1)$-dimensional Hausdorff measure of $F$. Consider the second order elliptic operator $\mathscr{L}^V:= -\sum_{l,k=1}^d \partial_l(c_{kl}\partial_k) + V$, where $c_{kl}=c_{lk}:\Om\to\R$ are real valued bounded smooth functions satisfying the uniform ellipticity condition
\begin{align*}
	C^{-1}|\xi|^2\leq\sum_{k,l=1}^dc_{k,l}(x)\xi_k\xi_l\leq C|\xi|^2,\quad \xi\in\R^d
\end{align*}
for some $C>0$
and $V:\Om\to\R$ is a real valued bounded non-negative measurable function. Assume further that the operator induced by $\mathscr{L}^V$ in $L^2(\Omega):=L^2(\Omega,dx)$ with Dirichlet boundary condition, say $L_D^V$, is injective. Let  $S_t^V$ be the semigroup of the Dirichlet-to-Neumann operator on $L^2(F):=L^2(F,\mu)$ related to $\mathscr{L}^V$ (see \cite{Ouhabaz} and references therein). Then $S_t^V$ maps $C(F)$ into itself and 
$S_t^V|_{C(F)}:C(F)\to C(F)$ has Feller property. 
\label{cor:Feller}
\end{coro}
\begin{proof}
Let $\calE^V$ be the symmetric regular Dirichlet form defined by
\[
\dom(\calE^V) = W^{1,2}(\Omega)\subset L^2(\overline\Omega):= L^2(\Omega,dx), \;
\calE^V(u,u) = \sum_{k,l=1}^d \int_\Omega c_{kl}|\nabla u|^2\,dx + \int_\Omega V u^2\,dx,
\]
and let $L^V$ be the related positive selfadjoint operator. We designate by $e^{-tL^V}$ the related semigroup. It is known that there is a unique standard process $\mathbb{X}^V$ with state space $\overline\Omega$ and symmetric semigroup $ {\bf P}_t^V$, associated to $\calE^V$ in the sense that there is a polar set $\mathscr{N}\subset\overline\Omega$ such that 
\[
e^{-tL^V} u(x) = \mathbb{E}_x\left[ e^{-\int_0^t V(\mathbb{X}_s)\,ds} \cdot u(\mathbb{X}_t) \right]={\bf{P}}_t^V u(x),\; x\in \overline{\Omega}\setminus\mathscr{N}, 
u\in\mathscr{B}_b(\overline\Omega).
\]
Here $\mathbb{X}=(\Omega,\mathbb{X}_t,\mathbb{P}_x)$ is the standard process associated to $\calE^V$ with $V=0$. 
For  $\alpha>0$ let $R_\alpha^V$ be the resolvent of ${\bf P}_t^V$. Then
\[
(L^V + \alpha)^{-1} u(x) = R_\alpha^V u(x),\; x\in \overline{\Omega}\setminus\mathscr{N},\; 
u\in\mathscr{B}_b(\overline\Omega).
\]
To ease notation we will drop the superscript $V$. By Sobolev inequality we conclude that $e^{-tL}, t>0$ and hence ${\bf P}_t$ has a symmetric absolutely continuous nonnegative kernel $p_t(\cdot,\cdot)$. Using \cite[Theorem 4.1]{Nittka} we derive that $p_t(\cdot,\cdot)$ can be continued into a H\"older continuous function on $\overline\Omega\times\overline\Omega$ which we still denote by $p_t(\cdot,\cdot)$, moreover $p_t(x,y)$ is continuous on $(0,\infty)\times\overline\Omega\times\overline\Omega$. Thus for any $\alpha>0$ both $(L + \alpha)^{-1}$ and $ R_\alpha$ have the same integral kernel, $G_\alpha(\cdot,\cdot)$. Moreover the function $G_\alpha(\cdot,y):\overline\Omega\to(0,\infty)$ is l.s.c for any fixed $y\in \overline\Omega$ by Fatou's Lemma.\\
Besides using  \cite[Theorem 6.10, p.~171]{Ouhabaz},  we get that $p_t$ satisfies the Gaussian upper bound
\[
p_t(x,y) \leq  \frac{c_1}{t^{d/2}} e^{c_2 t} e^{-\frac{|x-y|^2}{c_3 t}},\;\; t>0\;\; \text{and}\;\; x,y\in\overline\Omega,
\]
which in turn leads, for large $\alpha$,  to the fact that $G_\alpha(\cdot,\cdot)$ is continuous away from the diagonal and 
\[
G_\alpha (x,y) \leq c(\alpha) |x-y|^{2-d},\ x,y\in\overline\Omega.
\]	
Summing up  we get that $R_\alpha 1$ is continuous bounded and  $\mu$ is a $G_\alpha$-Kato measure for any $\alpha>0$, owing to the known estimate for the $(d-1)$-Hausdorff measure: $\mu(B(x,r)\cap F)\leq cr^{d-1}$ for all $x\in\R^d, 0<r\leq 1$ and all Euclidean open balls $B(x,r)$.\\
We claim that the fine support of $\mu$ is $F$. Indeed, the harmonic measure $H_F(x,\cdot)$ defined by 
\begin{align*}
	H_F(x,dy):=\mathbb{E}_{x}\left[e^{-\int_0^{\sigma_F}
		V(\mathbb{X}_s)ds}:\mathbb{X}_{\sigma_F}\in dy\right],\quad x\in\overline{\Omega}
\end{align*}
is absolutely continuous w.r.t $\mu(dy)$ by way of \cite[Theorem~3.1 a)]{JK}. We now prove that $F$ is the fine support of $\mu$ w.r.t $\mathbb{X}$. Let $u$ be a finely continuous nearly Borel measurable function on $\overline{\Omega}$ w.r.t the process $\mathbb{X}$. Assume that $u=0$ $\mu$-a.e. Then $H_F u=0$ on $\overline{\Omega}$, by the absolute continuity of $ H_F(x,dy) $ w.r.t $\mu(dy)$. Since $u=H_F u$ on $F$, we have $u=0$ on $F$. This means that $F$ is the fine support of $\mu$.\\    
Let $\bf{A}_t$ be the PCAF of $\mu$ w.r.t the process $\mathbb{X}$ and $\check{T_t}$ the time-changed semigroup of ${\bf P}_t$ by means of ${\bf A}_t$ considered on $\mathscr{B}_b(F)$ (see above). An appeal to Theorem \ref{Feller-restricted}-(1) implies that $\check{T_t}$  maps $C(F)$ into itself and $\check{T_t}:C(F)\to C(F)$ has Feller property.\\
Let us show that $\check{T_t}= S_t|_{C(F)}$. It is easy to check that $\check{T_t}$ defines a strongly continuous semigroup on $L^2(F)$. By \cite[Theorem 5.2.2, p.~177  and Corollary 2.2.15, p.~73]{CF} the regular Dirichlet form associated to the $L^2(F)$-semigroup $\check{T}_t$, which we designate by $\check\calE$ is given by
\[
\dom(\check\calE) =  \big\{u|_F\colon u\in L_{\rm loc}^2(\Omega),   \partial_i u\in L^2(\Omega), i=1,\cdots,d \big\},\   
\check\calE (u|_F,u|_F) = \calE(H_Fu, H_F u).
\]
The Dirichlet form $\check\calE$ is called the trace of $\calE$ on $F$ (see \cite[\S 5.2]{CF}).

Since $\Omega$ is bounded and Lipschitz, it is an extension domain so that  using Sobolev inequality we get $ W_{\rm ext}^{1,2} (\Omega) = W^{1,2} (\Omega)$. Using the known inequality 
\[
\int_F u^2\,d\mu \leq c\left(\int_\Om |\nabla u|^2\,dx + \int_\Om u^2\,dx\right),\ u\in W^{1,2}(\Omega),
\] 
we conclude that  $\dom(\check\calE) = H^{1/2}(F)$. 
Here $H^{1/2}(F)$ denotes the Sobolev space of order $1/2$ on $F$ defined by the trace operator from $W^{1,2}(\Omega)$ i.e. $H^{1/2}(F):=\{u\in L^2(F):\exists v\in W^{1,2}(\Omega)\text{ s.t. }
{\rm tr}(v)=u\}$.
Observing that for any $u\in W^{1,2} (\Omega)$ the function $H_F u$ is $L$-harmonic on $\Omega$  and  $u|_F = H_F u|_F$,  we conclude that  $\check\calE$ is exactly the Dirichlet form of the Dirichlet-to-Neumann operator as it is constructed in \cite{Ouhabaz} with associated semigroup $S_t$, which completes the proof.
\end{proof}

\begin{rk}
{\rm 
Corollary~\ref{cor:Feller} recovers the part which asserts that $S_t|_{C(F)}: C(F)\to C(F)$ defines a strongly continuous semigroup from theorem \cite[Theorem 1.1]{Ouhabaz}. However, with stronger assumption on the coefficients and weaker assumption on the domain. 
}
\end{rk}
\begin{theo}[Approximation]
Let $F_n,\, F$ be the fine supports of $\mu_n,\mu_\infty$, respectively. Assume that $F_n,\, F,\, n\in\N$ are closed. Let $\Pi_n : C_0(F)\to C_0(F_n),\;\, u\mapsto (P_F u)|_{F_n}$. Then for all $u\in C_0(F)$ it holds
\begin{enumerate}
\item[\rm(1)] 
\[
\lim_{n\to\infty}\sup_{x\in F_n} \left|\check{T}_t^n \Pi_n u(x) - \Pi_n \check{T}_t^\infty u(x) \right| = 0
\]
locally uniformly in $t\in [0,\infty)$.
\item[\rm(2)] For all $\alp>0$
\[
\lim_{n\to\infty}\sup_{x\in F_n} \left|V_\alp^n \Pi_n u(x) - \Pi_n V_\alp^\infty u(x) \right| = 0.
\]
\end{enumerate}
\label{Approximation}
\end{theo}		
Observe that for $u\in \mathscr{B}_b(F)$ we have $P_F u = P_F\tilde u$. We write simply $P_Fu$.
\begin{proof}
Let $u\in C_0(F)$ and $\tilde u\in C_0(X)$ an extension of $u$. Then for all $x\in F_n$ it holds
\begin{align*}
\check{T}_t^n \Pi_n u(x) - \Pi_n \check{T}_t^\infty u(x) &= \check{T}_t^n (P_F u|_{F_n})(x) - (P_F \check{T}_t^\infty u)|_{F_n} (x)\\
&= \check{T}_t^n P_F u(x) - P_F \check{T}_t^\infty u|_{F_n} (x)\\
& = \check{T}_t^n P_F \tilde u(x) - (P_F (\check{\bf P}_t^\infty \tilde u|_{F}))|_{F_n} (x)\\
& = (\check{\bf P}_t^n P_F \tilde u)|_{F_n}(x) - (P_F (\check{\bf P}_t^\infty \tilde u|_{F}))|_{F_n} (x)\\
&= (\check{\bf P}_t^n P_F \tilde u)|_{F_n}(x) - (P_F (\check{\bf P}_t^\infty \tilde u))|_{F_n} (x)\\
& = (\check{\bf P}_t^n P_F \tilde u)|_{F_n}(x) - \check{\bf P}_t^\infty \tilde u|_{F_n} (x).
\end{align*}
Using Theorem \ref{Convergence-Sg} we achieve 
\begin{align}
\sup_{x\in F_n} \left|\check{T}_t^n \Pi_n u(x) - \Pi_n \check{T}_t^\infty u(x) \right| 
\leq \|\check{\bf P}_t^n P_F \tilde u - \check{\bf P}_t^\infty \tilde u\|_\infty\to 0\quad \text{as}\quad n\to\infty
\end{align}
locally uniformly in $t\in [0,\infty)$.\\
The proof of the second assertion follows the same lines with the use of Theorem \ref{UniformConvergencePotOPerator}.
\end{proof}
Let us adopt the definition given in \cite[Definition 6.3, p.~146]{Engel-Nagel} concerning  mild solution of the heat equation. In the few following lines we omit the subscript $n$ for convenience.\\
Since $\check{T}_t$ has Feller property whenever $F$ is closed, then according to \cite[Proposition 6.4, p.~146]{Engel-Nagel} for any $v\in C_0(F)$ the function $\check{T}_t v$ is the unique mild solution of the heat equation
\begin{equation}
\begin{cases} 
		- \frac{\partial u(t)}{\partial t} = \check{L} u(t),\\
		u(0) = v.
\end{cases}
\end{equation}
Hence using Theorem \ref{Approximation}  we obtain: 
\begin{coro}
Under assumptions of Theorem \ref{Approximation}, let $v\in C_0(F)$ and  let $u_\infty$ be the mild solution of the heat equation
\begin{equation}
\begin{cases} 
		- \frac{\partial u_\infty(t)}{\partial t} = \check{L}^\infty u_\infty(t),\\
		u_\infty(0) = v
\end{cases}
\end{equation}
and for each $n\in\N$  let $u_n$ be the mild solution of the  heat equation
\begin{equation}
\begin{cases} 
		- \frac{\partial u_n(t)}{\partial t} = \check{L}^n u_n(t),\\
		u_n(0) =P_Fv|_{F_n}.
\end{cases}
\end{equation}
Then 
\begin{enumerate}
\item[\rm(1)]
\[
\lim_{n\to\infty}\|u_n(t) - P_F u_\infty(t) \|_{\infty,F_n} = 0
\]
locally uniformly in $t\in [0,\infty)$.
\item[\rm(2)] If moreover  $F_n\subset F$ for all $n$ then 
\[
\lim_{n\to\infty}\| u_n(t) -  u_\infty(t) \|_{\infty,F_n} = 0
\]
locally uniformly in $t\in [0,\infty)$.
\item[\rm(3)] Whereas if $F\subset F_n$ for all $n$ then 
\[
\lim_{n\to\infty}\| u_n(t) -  u_\infty(t) \|_{\infty,F} = 0
\]
locally uniformly in $t\in [0,\infty)$.
\end{enumerate}
Here  $\|\cdot\|_{\infty,F}$, $\|\cdot\|_{\infty,F_n}$ denotes the supremum respectivly  on $F$ and $F_n$.
\end{coro}
\begin{rk}
{\rm
The latter corollary recovers the result given by Ehnes and Hambly in \cite[Theorem 5.11]{Ehnes-Hambly}, where convergence of mild solutions and semigroups for some  Krein--Feller operators on the unit interval are considered.	

}
\end{rk}

We come now to convergence of finite distributions.
\begin{theo}\label{thm:finiteDimensionalConvergence}
For each $n\in\N_\infty$, let $\mathbb{E}^n$ be the expectation with respect to the process $\check{\mathbb{X}}^n$ with initial distribution $\mu_n$. Let $F_n,\,F$ be the fine support of $\mu_n,\,\mu_{\infty}$, respectively.  
Assume that $F_n,\,F$ are closed and either $F_n\subset F$ for all $n$ or $\lim_{n\to\infty} \|P_{F_n} u - P_F u\|_\infty =0$ for all $u\in C_0(X)$. Then the finite time distributions of $\check{\mathbb{X}}^n$ converge vaguely to the finite time distribution of $\check{\mathbb{X}}^\infty$, i.e., for every integer $k$, every $0<t_1<\cdots<t_k<\infty$ and every continuous function with compact support $u:X^{k+1}\to\R$ it holds
\[
\lim_{n\to\infty}\mathbb{E}^n \left[u\left( \check{\mathbb{X}}^n(0), \check{\mathbb{X}}^n(t_1),\cdots,\check{\mathbb{X}}^n(t_k)\right)\right]
	= \mathbb{E}^\infty \left[u\left( \check{\mathbb{X}}^\infty(0), \check{\mathbb{X}}^\infty(t_1),\cdots,\check{\mathbb{X}}^\infty(t_k)\right)\right].
\]
\end{theo}
\begin{proof}
The proof is inspired from \cite{HZ}.\\
By Stone--Weierstrass theorem it suffices to consider functions $u$ of the type 
\[u(x_0,x_1,\cdots,x_k) = u_0(x_0)u_1(x_1)\cdots u_k(x_k)
\quad\text{ with }\quad u_i\in C_c(X)\quad\text{ for\; each }\quad i.
\] 
For such type of  functions  $u$ it holds
\begin{align*}
		\mathbb{E}^n& \left[u\left( \check{\mathbb{X}}^n(0), \check{\mathbb{X}}^n(t_1),\cdots,\check{\mathbb{X}}^n(t_1)\right)\right]
		 - \mathbb{E}^\infty \left[u\left( \check{\mathbb{X}}^\infty(0), \check{\mathbb{X}}^\infty(t_1),\cdots,\check{\mathbb{X}}^\infty(t_k)\right)\right]\\
		& =\int_X u_0 \check{{\bf P}}_{t_1}^n( u_1  \check{{\bf P}}_{t_2 - t_1}^n( u_2\cdots   u_{k-1} \check{{\bf P}}_{t_{k} - t_{k-1}}^n u_k ) )\,d\mu_n\\
		&\quad- \int_X u_0 \check{{\bf P}}_{t_1}^\infty( u_1  \check{{\bf P}}_{t_2 - t_1}^\infty( u_2\cdots   u_{k-1} \check{{\bf P}}_{t_{k} - t_{k-1}}^\infty u_k ) )\,d\mu_\infty\\
		& = \int_X u_0 \check{{\bf P}}_{t_1}^n( u_1  \check{{\bf P}}_{t_2 - t_1}^n( u_2\cdots   u_{k-1} \check{{\bf P}}_{t_{k} - t_{k-1}}^n u_k ) )\,d\mu_n\\
		&\quad- \int_X u_0 \check{{\bf P}}_{t_1}^\infty( u_1  \check{{\bf P}}_{t_2 - t_1}^\infty( u_2\cdots   u_{k-1} \check{{\bf P}}_{t_{k} - t_{k-1}}^\infty u_k ) )\,d\mu_n\\
		& +  \int_X u_0 \check{{\bf P}}_{t_1}^\infty( u_1  \check{{\bf P}}_{t_2 - t_1}^\infty( u_2\cdots   u_{k-1} \check{{\bf P}}_{t_{k} - t_{k-1}}^\infty u_k ) )\,d\mu_n\\
		&\quad - \int_X u_0 \check{{\bf P}}_{t_1}^\infty( u_1  \check{{\bf P}}_{t_2 - t_1}^\infty( u_2\cdots   u_{k-1} \check{{\bf P}}_{t_{k} - t_{k-1}}^\infty u_k ) )\,d\mu_\infty.
\end{align*}
The first difference tends to zero as $n\to\infty$ by uniform convergence of $\check{{\bf P}}_{t}^n$ towards $\check{{\bf P}}_{t}^\infty$ on $C_0(X)$, the fact that $u_0$ has compact support and the vague convergence of $\mu_n$ towards $\mu_\infty$. Also the second difference tends to zero as  $n\to\infty$, owing to vague convergence of $\mu_n$ and the fact that $u_0$ has compact support.
\end{proof}

Finally we give the weak convergence of time-changed processes. 
\begin{theo}\label{thm:weakconvergence}
For each $n\in\N_\infty$, let $\mathbb{P}^n$ be the probability law on 
the path space $D_{F_{\partial}}[0,\infty)$ of right continuous paths having left-hand limit in $F_{\partial}$ 
with respect to the process $\check{\mathbb{X}}^n$ with initial distribution $\mu_n$. Let $F_n,\,F$ be the fine support of $\mu_n,\,\mu_{\infty}$, respectively.  
Assume that $F_n,\,F$ are closed, $F_n\subset F$ for all $n\in\mathbb{N}$. Then $\check{\mathbb{X}}^n$ converges weakly $\check{\mathbb{X}}^\infty$.
\end{theo}
\begin{proof}
For $n\in\N$, let $S_t^n$ be the semigroup defined by
\[
S_t^n : C_0(F)\to C_0(F),\quad S_t^n u(x) = \mathbb{E}_x\big ( u|_{F_n}(\mathbb{X}_{\tau_t^n}) \big).
\]	
By Theorem~\ref{thm:Feller}, we conclude that $S_t^n$ has Feller property on each $F_n$. Note here that $S_t^n$ preserves $C_0(F)$ itself and  also has  Feller property on $F$ in view of Theorem~\ref{Substitute-Feller}-(3) and $S_0^n=P_{F_n}$ by replacing $X$ with $F$. So we can deduce the weak convergence by applying \cite[2.5 Theorem, p.~167]{EK}. We outline the proof: 
Let $\mathbb{Y}^n$ be the standard process associated to $S_t^n$. Then the sample paths of $\mathbb{Y}^n$ coincide with those of $\check{\mathbb{X}}^n$ with probability $1$ 
under $\mathbb{P}^x$ for $x\in F_n$ and hence they have the same law. 
On the other hand by Proposition~\ref{prop:ConvergenceFeller}, for $u\in C_0(F)$, we have the uniform convergence
\begin{align}                                                                                        
\lim_{n\to\infty}\|S_t^n u-\check{\bf P}_t^\infty u\|_{\infty,F}=0\label{eq:uniformconvergence}
\end{align}
for all $t\geq0$ and the convergence is uniform in $t$ on compact intervals of $[0,\infty)$. Consequently one can apply \cite[2.5 Theorem, p.~167]{EK} to obtain weak convergence of the processes.

\end{proof}

\noindent{\bf Declarations.} The authors have no conflicts of interest to declare that are relevant to the content of this article. \\ 

\noindent{\bf Ethical approval.} Not applicable.\\
\quad\\
{\bf Funding.} No funding for the first named author, the second named author is supported in part by JSPS Grant-in-Aid for Scientific Research (S) (No. 22H04942).\\
\quad\\
{\bf Availability of data and materials.}  Data sharing is not applicable to this article as no datasets were generated or analyzed during the current study.
\\
\quad
\\
\noindent{\bf Acknowledgement.} The authors would like to thank Professors 
Yuichi Shiozawa and Kouhei Matsuura for suggesting 
to utilize \cite[2.5 Theorem, p.~167]{EK} for the weak convergence of Feller processes. They also thank the anonymous referee, whose comments improve the quality of this paper so much.

\printbibliography


\end{document}